\documentclass[a4paper]{article}

\usepackage[all]{xy}
\usepackage{epsf} 
\usepackage{amsfonts}
\usepackage{amssymb}
\newtheorem{st}      {Theorem}    [section]
\newtheorem{prop}[st]{Proposition}

\newtheorem{exam}[st]{Example}   
\newtheorem{num}[st]{}

\newcommand{\ii}{\ensuremath{\mathcal{I}}}
\newcommand{\thp}{\ensuremath{\theta}}
\newcommand{\Ga}{\ensuremath{\Gamma}}
\newcommand{\ga}{\ensuremath{\gamma}}
\newcommand{\al}{\ensuremath{\alpha}}
\newcommand{\be}{\ensuremath{\beta}}
\newcommand{\f}{\ensuremath{\varphi}}
\newcommand{\rmap}{\longrightarrow}
\newcommand{\lmap}{\longleftarrow}
\newcommand{\Ka}{\ensuremath{\mathcal{K}}}
\newcommand{\Ha}{\ensuremath{\mathcal{H}}}
\newcommand{\M}{\ensuremath{\mathcal{M}}}
\newcommand{\N}{\ensuremath{\mathcal{N}}}

\newcommand{\U}{\ensuremath{\mathcal{U}}}
\newcommand{\OO}{\ensuremath{\mathcal{O}}}

\newcommand{\la}{\ensuremath{\Lambda}}

\newcommand{\A}{\ensuremath{\mathcal{A}}}

\newcommand{\el}{\ensuremath{\mathcal{L}}}
\newcommand{\F}{\ensuremath{\mathcal{F}}}
\newcommand{\es}{\ensuremath{\mathcal{S}}}
\newcommand{\B}{\ensuremath{\mathcal{B}}}
\newcommand{\C}{\ensuremath{\mathcal{C}}}
\newcommand{\G}{\ensuremath{\mathcal{G}}}
\newcommand{\g}{\ensuremath{\Gamma}}
\newcommand{\gc}{\ensuremath{\Gamma_c}}
\newcommand{\ps}{{\raise 1pt\hbox{\tiny (}}}
\newcommand{\pss}{{\raise 1pt\hbox{\tiny [}}}
\newcommand{\pdd}{{\raise 1pt\hbox{\tiny ]}}}
\newcommand{\pd}{{\raise 1pt\hbox{\tiny )}}}
\newcommand{\bs}{{\raise 1pt\hbox{\tiny [}}}
\newcommand{\bd}{{\raise 1pt\hbox{\tiny ]}}}
\newcommand{\nG}[1]{\ensuremath{\G^{\ps #1\pd}}}

\newcommand{\nXG}[1]{\ensuremath{(X \cross \G)^{\ps #1\pd}}}
\def\cross{\mathinner{\mathrel{\raise0.8pt\hbox{$\scriptstyle>$}}
                 \joinrel\mathrel\triangleleft}}
\newcommand{\nH}[1]{\ensuremath{\Ha^{\ps #1\pd}}}
\newcommand{\nB}[1]{\ensuremath{B^{\ps #1\pd}}}

\newcommand{\gth}{\ensuremath{\G_{\ps \theta \pd}}}

\newcommand{\conv}{\ensuremath{\C_{c}^{\infty}(\G)}}
\newcommand{\Z}{\ensuremath{\mathcal{Z}}}
\newcommand{\nZ}[1]{\ensuremath{\Z^{\ps #1\pd}}}

\def\dt{{\raise 0.5pt\hbox{$\scriptscriptstyle\bullet$}}}
\def\compose{{\raise 1pt\hbox{$\scriptscriptstyle\circ$}}}
\def\dcross{{\raise 0.5pt\hbox{$\scriptscriptstyle\boxtime$}}}

\begin{document}

\title {Cyclic Cohomology of \'Etale Groupoids; The General Case}
\author {Marius Crainic}
\date { Mail address:Utrecht University, Department of Mathematics, P.O.Box:80.010, Budapestlaan 6, 3508 TA Utrecht, The Netherlands
\\ e-mail address: crainic@math.ruu.nl}
\pagestyle{myheadings}
\maketitle
\begin{abstract}
We give a general method for computing the cyclic cohomology of crossed products by \'etale groupoids, extending the Feigin-Tsygan-Nistor spectral sequences. In particular we extend the computations performed by Connes, Brylinski, Burghelea and Nistor for the convolution algebra $\conv$ of an \'etale groupoid, removing the \emph{Hausdorffness} condition and including the computation of \emph{hyperbolic} components. Examples like group actions on manifolds and foliations are considered.
\par {\bf Keywords}: cyclic cohomology, groupoids, crossed products, duality, foliations.
\end{abstract}

\tableofcontents
\newpage

\section{Introduction}

\ \ \ \ In the general picture of non-commutative geometry, cyclic homology plays the role of compactly supported de Rham cohomology; it was introduced by A. Connes as the target of the Chern character. The dual theory is cyclic cohomology, which plays the role of closed de Rham homology. The pairing between these two is an important tool in performing numerical computations of K-theory classes (indices).
\par Often the non-commutative space we have to deal with is an orbit space of an \'etale groupoid; in particular, any \'etale groupoid can be viewed as such a non-commutative space. This fits in with Grothendieck's idea of what a ``generalized space'' is (\cite{SGA,Mo1}), and includes examples like  leaf spaces of  foliations, orbit spaces of group actions on manifolds, orbifolds. To say what the groups $HP_{*}(\conv)$ look like is an important step in solving index problems and in understanding the connection between the topology and the analysis of ``leaf spaces'' (here we have in mind in particular the Baum-Connes assembly map \cite{BaumConnes}). 
\par The computation of $HP_{*}(\conv)$ was started by Connes for the case where $\G= M$ is a manifold (\cite{Co2}), Burghelea for the case where $\G=G$ is a group (\cite{Bu}) and by Feigin, Tsygan and Nistor for crossed products by groups (\cite{FeTs}, \cite{Ni1}). The general strategy is to decompose these homology groups as direct sums of localized homologies; there are two different kinds of components which behave differently. Following the terminology introduced in \cite{Bu}, these are called elliptic and hyperbolic components. Usually the hyperbolic ones are more difficult to compute and involve in a deeper way the combinatorics of the groupoid.
\par In the general setting of smooth \'etale groupoids the results were partially extended by Brylinski and Nistor (\cite{BrNi}): for a \emph{Hausdorff} \'etale groupoid \G, the localized homologies $HP_{*}(\conv)_{\OO}$ are  defined for any invariant open space $\OO$ of loops; the \emph{elliptic} components are computed in terms of double complexes; in particular the localization at units is related to the homology of the classifying space $B\G$. There are some important questions left:
\par 1) Compute the hyperbolic components;
\par 2) Remove the Hausdorffness condition (simple examples coming from foliations are non-Hausdorff);
\par 3) Find a way to book-keep the computations; in particular give a more conceptual proof and a more conceptual meaning of the results.
\par The initial intention of the author (after studying the computations made in \cite{BrNi}) was to answer the last question; this led to the definition of a homology theory for \'etale groupoids (which could be called cohomology with compact support as well) developed in \cite{CrMo}. This homology turned out to be extremely helpful in answering all the questions above, among others. As it is very flexible, we expect it to be useful also in understanding the Baum-Connes assembly map.
\par The approach and the results of this paper owe a great deal to the previous work of several authors, especially Brylinski (\cite{BrNi}), Burghelea (\cite{Bu}) and Nistor (\cite{BrNi}, \cite{Ni1}). The computation we give for the localization at units is, beyond the formalism, the same as the one given in \cite{BrNi}; to the same paper we owe the important idea of reduction to loops (proposition \ref{redloops}). The method for computing the other localizations are inspired by the initial work of Burghelea \cite{Bu}; the difficulty is that the topological arguments used in that paper do not work in this generality any more. An older idea (\cite{Mo1}) that working with classifying toposes (i.e. sheaves) might be easier than working with classifying spaces (and this was pointed out, for the first time in our context, in the same paper \cite{BrNi}) becomes essential for us. With this in mind, our job is to replace the classifying spaces used by Burghelea by suitable sheaves on \'etale categories and the topological arguments by a suitable algebraic-topological formalism (long exact sequences and spectral sequences for homology of \'etale groupoids).
\par The paper is organized as follows. In section 2 we review the basic definitions and properties of sheaves on \'etale groupoids and their homology and show how the Hausdorffness assumption can be dropped. It is important to point out that our definition (see \ref{our}) of compactly supported forms on non-Hausdorff manifolds (which was first introduced in \cite{CrMo}) is related to, but not the same as the one given by Connes (section 6 in \cite{Connes}). Ours has basic properties, like the existence of a de Rham differential, which are not shared by Connes' (as remarked in the introduction of \cite{BrNi}); it is also the right object for extending Poincar\'e-duality to non-Hausdorff manifolds (\cite{CrMo}). For this reasons we expect it to be useful also in other problems which deal with foliations with non-Hausdorff graph. 
\par In section 3 we introduce the homology of groupoids with coefficients in cyclic sheaves (subsection 3.2); more generally, given a cyclic groupoid (i.e. a groupoid with an action of $\mathbb{Z}$ on it, see \ref{cycat}) we consider twisted cyclic sheaves (for which the usual identity $t^{\ps n+1\pd}=1$
is replaced by  $t^{\ps n+1\pd}=$the action of the generator $1 \in \mathbb{Z}$, see \ref{twist}). In subsection 3.3 we prove the main technical results  concerning these homologies; at the end we derive as a simple consequence the Eilenberg-Zilber-type spectral sequence for cyclic objects which is one of the main results in \cite{GeJo}. The older approaches to cyclic homology of crossed products by (discrete) groups can not be directly extended to the setting of \'etale groupoids; we show in subsection 3.4 how cyclic groupoids can be used to overcome this problem. In particular we extend the Feigin-Tsygan-Nistor spectral sequences (\cite{FeTs,Ni1}) and Nistor's description of the $S$ boundary (\cite{Ni1}). See theorem \ref{1elliptic}, \ref{1hyperbolic}\ .
\par In section 4 we come down to earth with more concrete applications; here is a list of them:
\par 4) For smooth \'etale groupoids we extend the old results of Burghelea proving that the elliptic components $HP_{*}(\conv)_{\OO}$ are computed by the homology of the normalizer $\N_{\OO}$ of $\OO$ (see theorems \ref{ultima}, \ref{last});
\par 5) For hyperbolic components $HC_{*}(\conv)_{\OO}$, we describe a $H^{*}(\N_{\OO})$-module structure which identifies $S$ in the $SBI$-sequence with the product by an element $e_{\OO} \in H^{2}(\N_{\OO})$; in particular we get a vanishing condition for $HP_{*}(\conv)_{\OO}$ (which extends a similar result of Burghelea \cite{Bu} and Nistor \cite{Ni1}). For stable $\OO$'s we also give a more concrete description of $HC_{*}(\conv)_{\OO}$. See theorem \ref{lastt};
\par 6) In subsection 4.3 we show how the methods apply to \emph{co}homology. In particular we get that the pairing between $HP_{*}$ and $HP^{*}$ is a Poincar\'e-duality pairing, so it is highly non-trivial. See \ref{cohelliptic}, \ref{of}\ .
\par 7) For group-actions on manifolds we get the old results for the elliptic components (\cite{BaCo}, \cite{BrNi}), and a new description of the hyperbolic ones (see Corollary \ref{off})\ ;
\par 8) For foliations we prove that the cyclic homology is a well defined invariant of the leaf space of the foliation, in the sense that the process of reducing to the setting of smooth \'etale groupoids does not depend on the choice of the complete transversal (see theorem \ref{folinv}). We also give some examples. 
\par {\bf Acknowledgments}: First of all I would like to thank my supervisor, Prof. I.Moerdijk for many helpful discussions. Then I would like to thank ``Moerdijk's geometry team'' at Utrecht University for a nice working environment and useful coffee breaks; among them are Utpal Sarkar(with whom I studied the paper \cite{BrNi} and formulated Proposition \ref{redloops}), Bella Baci (who gave me many interesting suggestions) and Japie Vermeulen (who tried to teach me some topos theory). I would like to thank also Prof. Victor Nistor for his preparedness to help and for his patience in answering all my e-mails.
\par Finally I would like to thank the NWO (The Netherlands Organization for Scientific Research) which is financing my PhD research.


\section{Homology and Cohomology of Sheaves on \'Etale Groupoids}


\subsection{\'Etale Groupoids}

\begin{num}
\emph{
{\bf Groupoids} : A groupoid is a small category \G\ in which every arrow is invertible;
so it is given by a set \nG{0}\ of ``objects'', a set \nG{1}\ of ``arrows'' and maps $s,t: \nG{1} \longrightarrow \nG{0}$ for source and target, $u: \nG{0} \longrightarrow \nG{1}$ for units, $i:\nG{1} \rmap \nG{1}$ for inverse and 
$m:\nG{2}=\{ (g,h)\in\nG{1}\times\nG{1}\ : s(g)=t(h)\} \longrightarrow\nG{1}$ for
 multiplication. A topological groupoid is a groupoid \G\ with topologies on 
\nG{0}\ and \nG{1}\ such that all the structure maps $s,t,u,i,m$ are continuous. It is called \'etale if $s$ is an \'etale map (i.e. a local homeomorphism). In this case all the other structure maps are \'etale; in particular \nG{0}\ is open in \nG{1}\ (see also $2.7$ in \cite{BrNi}). A smooth groupoid is a groupoid \G\ with differentiable structures on \nG{0}\ and \nG{1}\ such that \nG{0}\ is Hausdorff and all the structure maps are smooth. It is called smooth \'etale if $s$ is a local diffeomorphism. We shall denote $m \ps g,h \pd =gh ,\ i \ps g \pd =g ^{-1}$
and $g:c\longrightarrow d$ to denote $s \ps g \pd =c ,\ t \ps g \pd =d.$ The space of arrows is also denoted by \G.
}
\end{num}

\begin{num}
\label{exgr}
\emph{
{\bf Examples}
\begin{enumerate}
\item  Any space $X$ can be viewed as an \'etale groupoid with $\nG{0}=\nG{1}=X.$
\item Any topological group G can be viewed as a topological groupoid with one object, G as the space of arrows and with the multiplication of G. It is \'etale if and only if $G$ is discrete.
\item As a mixture of the previous examples, if $G$ is a group acting on a space $X$, the cross-product groupoid $X \cross G$ is defined by {\ensuremath{(X \cross G)^{\ps 0\pd}}} $=X$, {\ensuremath{(X \cross G)^{\ps 1\pd}}} $=X \times G$, 
$s \ps x,g \pd =xg , t \ps x,g \pd =x , u \ps x \pd = \ps x,1 \pd , m( \ps x,g \pd ,\ps y,h \pd )= \ps x,gh \pd, i\ps x, g \pd=\ps xg, g^{-1} \pd$. It is a good replacement for the orbit space $X/G$ (see \cite{Co3}).
\item Many examples of groupoids arise in foliation theory: Haefliger's groupoid $\g ^q$, or the holonomy groupoid $Hol(M,\F)$ of a foliated manifold $(M,\F)$. The latter is \'etale if one reduces the space of objects to a complete transversal, and is a good replacement for the leaf space of the foliation. See e.g. \cite{Ha2} , \cite{Wi} , \cite{Connes} .
\item Orbifolds can be modelled by \'etale groupoids ; they correspond to the \'etale groupoids \G\ with the property that $ (s,t):\G \longrightarrow \nG{0} \times \nG{0} $ is a proper map (see \cite{MoPr}) .
\end{enumerate}
}
\end{num}

\begin{num}
\label{actions}
\emph{
{\bf Actions}: Let \G\ be an \'etale groupoid. A right action of \G\ on the space $X$ consists of two continuous maps $\pi : X \longrightarrow \nG{0}$ (the moment map), $m:X \times_{\nG{0}} \G =\{ \ps x,g \pd \in X \times \G : \pi \ps x \pd =t \ps g \pd \} \longrightarrow X$
(the action)  such that, denoting $m \ps x,g \pd =xg$:\[\ps xg \pd h=x \ps gh \pd , x1=x,\pi \ps xg \pd =s \ps g \pd .\]
}
\end{num}
\ \ \ We shall call $X$  a right $\G$ space with  moment map $\pi$. The associated groupoid for this action, denoted $X \cross \G $, is defined as a generalization of $1.2.3 $: $ \nXG{0} = X , \nXG{1} = X \times_{\nG{0}} \G , s \ps x,g \pd =xg, t \ps x,g  \pd =x, u\ps x \pd  = \ps x,1 \pd , m(\ps x,g \pd,\ps y,h \pd ) =\ps x,gh \pd, i\ps x, g \pd=\ps xg, g^{-1} \pd$.
\par There is an obvious similar notion of left \G-space. Unless specified, all \G-spaces will have the action from the right.


\begin{num}
\label{bun}
\emph{
{\bf Bundles}: A (right) \G-bundle over the space B consists of a \G-space $E$ and a continuous map $p:E \longrightarrow B$ which is \G-invariant (i.e. $p \ps xg \pd =p \ps x \pd $). It is called principal if $p$ is an open surjection and 
$E \times_{\nG{0}}\G \longrightarrow E \times_{B} E , \ps e,g \pd  \mapsto \ps e,eg \pd$ is a homeomorphism.
}\end{num}

\begin{num}
\label{morph}
\emph{
{\bf Morphisms of Groupoids} (\cite{Haf},\cite{Mo1}): Let \G\ and \Ha\  be two groupoids. A morphism $ P:\G \longrightarrow \Ha $ from \G\ to \Ha\  (or Hilsum-Skandalis map cf. \cite{Mr}) consists of a space $P$, continuous maps (source and target): $s_{P}: P\longrightarrow \nG{0} , t_{P}: P \longrightarrow \nH{0}$, a left action of \G\ on $P$ with the moment map $s_{P}$, a right action of \Ha\ on $P$ with the moment map $t_{P}$, such that: \begin{enumerate}
\item $s_{P}$ is \Ha-invariant, $t_{P}$ is \G-invariant;
\item the actions of \G\ and \Ha\ on $P$ are compatible: $\ps gp \pd h=g \ps ph \pd $;
\item $s_{P}: P\longrightarrow \nG{0}$, as an \Ha-bundle with the moment map $t_{P}$, is principal.
\end{enumerate}
\ \ \ A nice intuitive motivation of this definition is that $P$ can be viewed as a continuous map between the orbit spaces of $\G$ and $\Ha$, described by its graph (see II.8.$\ga$ in \cite{Connes} or \cite{Ha3}). A nice theoretical motivation is that these morphisms are exactly the topos-theoretic morphisms between the orbit spaces of $\G$ and $\Ha$ viewed as toposes (i.e. between the classifying toposes of $\G$ and $\Ha$; see \cite{Mo1} for the precise statements and descriptions). The composition of two morphisms $ P:\G \longrightarrow \Ha $, $Q:\Ha \longrightarrow \Ka $ is defined by dividing out $P \times_{\nH{0}} Q$ by the action of \Ha: $ \ps p,q \pd h= \ps ph,h^{-1}q \ \pd $, and taking the obvious actions of \G\ and \Ha\ . We get in this way the category of groupoids and its full subcategory of \'etale groupoids.
}\end{num}

\begin{num}
\label{exf}
\emph{
{\bf Example}: Any continuous functor $\f :\G \rmap \Ha$ can be viewed as a morphism by taking $P_{\f}=\nG{0} \times_{\nH{0}} \Ha =\{ \ps c,h \pd \in \nG{0} \times \Ha : \f \ps c \pd =t \ps h \pd \},\ s_{P} \ps c,h \pd =c, \ t_{P} \ps c,h \pd =s \ps h \pd$ and the obvious actions.
}\end{num}

\begin{num}
\emph{
{\bf Morita Equivalences}: Two groupoids \G\ and \Ha\ are called Morita equivalent if they are isomorphic in the category of groupoids (as defined in \ref{morph}). An isomorphism $ P:\G \longrightarrow \Ha $ is called a Morita equivalence (cf. \cite{Mo1},\cite{Ha3}).
}
\end{num}

\begin{num}
\label{eseq}
\emph{
{\bf Examples}:
\begin{enumerate}
\item Recall (\cite{Mo1}) that an essential equivalence is a continuous functor $\f :\G \rmap \Ha$ with the property that $P_{\f}=\nG{0} \times_{\nH{0}} \Ha =\{ \ps c,h \pd \in \nG{0} \times \Ha : \f \ps c \pd=t \ps h \pd \} \rmap \nH{0} , \ \ps c,h \pd \mapsto s \ps h \pd$ is an open surjection and the diagram:
\[ \xymatrix {
\G \ar[r] ^-{\f} \ar[d] _-{\ps s,t \pd} & \Ha \ar[d] ^-{\ps s,t \pd} \\
 \nG{0} \times \nG{0} \ar[r] ^-{\f \times \f} & \ \nH{0} \times \nH{0} } \]
is a pull-back of topological spaces. It is easily seen that in this case \f\ induces a Morita equivalence $P_{\f}:\G \longrightarrow \Ha $ (see \ref{exf}).
   In fact we can prove (see 2.3 in \cite{Mo1} or \cite{Mr}) that \G\ and \Ha\ are Morita equivalent if and only if there is a groupoid \Ka\ and essential equivalences $\G \lmap \Ka \rmap \Ha$. 
\item If $(M,\F)$ is a foliated manifold, $T \rmap M$ is a complete transversal (\cite{Connes,Co4,HiSk}), there is an obvious functor $Hol_{T}(M,\F) \rmap Hol(M,\F)$ from the holonomy groupoid restricted to $T$ to the holonomy groupoid. It is a standard fact that this is a Morita equivalence.
\item If $E \rmap B$ is a principal $G$-bundle (where $G$ is a topological group), then the obvious projection $E \cross G \rmap B$ (see  examples 1, 3 in \ref{exgr}) is a Morita equivalence.
\item If \G\ is an \'etale groupoid, $\U=\{U_{\al}:\al \in I \}$ is a covering of \nG{0}, the groupoid $\G_{\U}$ is defined by \[\nG{0}_{\U}=\bigcup_{\al \in I} U_{\al} \times \{ \al \} \ , \ \nG{1}_{\U}=\bigcup_{\al , \be \in I} (s ^{-1} \ps U_{\al} \pd \cap t ^{-1} \ps U_{\be} \pd ) \times \{ \al , \be \}, \]
 $s \ps x,\al ,\be \pd =(s \ps x \pd , \al ) ,  t \ps x,\al ,\be \pd =(t \ps x \pd ,\be )$. The obvious functor $\G_{\U} \rmap \G$ is a Morita equivalence (\cite{Ha3}).
\end{enumerate} 
}
\end{num}

\begin{num}
\label{comma}
\emph{
{\bf Comma Groupoids}: If $ \f:\G \longrightarrow \Ha $ is a continuous functor between groupoids, $d \in \nH{0}$, the comma groupoid $d/\f$ is defined as follows. It has as objects pairs $ \ps h, c \pd \in \Ha \times \nG{0}$ with $s \ps h \pd =d$, $t\ps h \pd=\f \ps c \pd$ (i.e. the space of objects is $\nG{0} \times_{\nH{0}} \nH{1}$), and as morphisms from $\ps h,c \pd$ to $\ps h\ ',c\ ' \pd$ those $g : c \rmap c\ ' $ in \G\ with $\f \ps g \pd h = h\ ' $ (i.e. the space of morphisms is $\G \times_{\nG{0}} \Ha$). We have a commutative diagram:
\[ \xymatrix {
 d/\f \ar[r]^-{\omega_{d}} \ar[d] & \G \ar[d]^-{\f} \\
 \dt \ar[r] \ar[r]^-{d} & \Ha } \]
where $\omega_{d}$ is the continuous functor which send an object $\ps h, c\pd$ to $c$ and a morphism $g$, from $\ps h,c \pd$ to $\ps h\ ',c\ ' \pd$, to $g$. The comma category can be viewed as the fiber of $\f$ above $d$.
}
\end{num}


\subsection{\gc\ in the non-Hausdorff case}

\ \ \ It is well known in sheaf theory that usual notions concerning compactness do not behave well on non-Hausdorff spaces. In the sequel, the notions of c-softness, compactly supported cohomology, the compactly supported sections functor \gc\ or, more generally, the functor $f_{\, !} :Sh(X) \rmap Sh(Y)$ induced by a continuous map $f :X \rmap Y$ (see \cite{Go,Iv}) will have the usual meaning only on Hausdorff spaces. For locally Hausdorff spaces (i.e. spaces which have a Hausdorff open covering) a good extension of these notions is developed in \cite{CrMo}. In this subsection we briefly recall the main definitions and properties.

\begin{num}
\label{ass}
{\bf Assumption}:Throughout this paper, all spaces are assumed to be locally Hausdorff, locally compact, of finite cohomological dimension ;  all sheaves are assumed to be sheaves of complex (or real) vector spaces.
\end{num}
By ``of finite cohomological dimension'' we mean there is an integer $n$ such that every point has a Hausdorff neighborhood of cohomological dimension $\leq n$. The smallest $n$ with this property is called the cohomological dimension of the given space.

\begin{num}
\emph{
{\bf c-Softness, \gc}: Let $X$ be a space. A sheaf $\A \in Sh(X)$ is called c-soft if there is a Hausdorff open covering \U\ of $X$  such that $\A |_{U} \in Sh(U)$ is c-soft for all $ U \in \U$. In this case define $\gc(U, \A)$ as the image of the map:
       \[ \bigoplus_{U} \gc (U;\A) \rmap \Ga (X_{dis} ;\A)   \ \ \ \ \ ( U \in \U) ,\]
}
\end{num}
where $\Ga (X_{dis};\A)=\{ u:X \rmap \bigsqcup_{x \in X} \A_{x} : u \ps x \pd \in \A_{x}  \ \ ,\forall  \ x \in \A_{x} \}$ \ ($X_{dis}$ is $X$, considered with the discrete topology)
and $\gc (U;\A) \rmap \Ga (X_{dis} ;\A) $, $s \rmap \overline{s}$ is given by:
      \[ \overline{s}(x)=germ_{x}(s)  \ \  {\rm for} \ x \in U \ , \ \ \ {\rm and} \ 0 \ {\rm otherwise}. \]
\ \ \ The basic property which enables us to extend the usual results from the Hausdorff case is the Mayer-Vietoris sequence (like in \cite{BoTu}, pp.139,186): for any open covering \U\ of $X$, there is a long exact sequence:
\[  . \ . \ . \rmap \bigoplus_{U,V} \gc(U \cap V ;\A) \rmap \bigoplus_{U} \gc(U;\A) \rmap \gc(X ; \A) \rmap 0 \ \  .\]
\ \ \ The compactly supported cohomology functor, $H^{*}_{c}(X,-)$ is defined by using c-soft resolutions (as in \cite{Go,Iv}); in particular we define $\gc(X,\es)=H^{0}_{c}(X,\es)$ for any sheaf $\es \in Sh(X)$.
\par 
For any continuous map $f:X \rmap Y$ and any c-soft sheaf $\A \in Sh(X)$, there is an unique c-soft sheaf $f_{\, !}\A \in Sh(Y)$ such that we have the ("usual") isomorphisms $\gc(V;f_{\, !}\A) \simeq \gc(f^{-1}(V);\A) $ for all $\ V \subset\ Y$ open, functorial with respect to the inclusion of opens in $Y$. Using c-soft resolutions in $Sh(X)$ one can define $f_{\, !}(\es) \in Sh(Y)$, or more generally $R^{*}f_{\, !}(\es) \in Sh(Y)$ for any $\es \in Sh(X)$ (see also \cite{Go,Iv}).

\begin{num}
\label{basicsgc}
\emph{
{\bf Basic properties of \gc,\ $f_{ \, !}$}: The usual properties of \gc,\ $f_{\, !}$ in the Hausdorff case (\cite{Go,Iv}) easily extend to our case by using the  Mayer-Vietoris sequence:
\begin{enumerate}
\item $\gc(X,-)$ preserves quasi-isomorphisms between bounded below cochain complexes of c-soft sheaves. In fact, under our assumption (that $X$ has finite cohomological dimension) this is true also for unbounded complexes.
\item If $A \subset X$ is closed, $\A \in Sh(X)$ is c-soft, then the obvious restriction $\gc(X,\A) \rmap \gc(A,\A)$ induces an isomorphism $\gc(A,\A) \simeq \gc(X,\A)/ \{u \in \gc(X,\A): u|_{A}=0 \}$. If $A=\theta^{-1} \ps 0 \pd$ for some continuous map $\theta :X \rmap \mathbb{R}, u \in \gc(X,\A)$ then:
\[ u|_{A}=0 \Longleftrightarrow   \ \exists \epsilon > 0  :u|_{\theta^{-1}(\ps -\epsilon , \epsilon \pd)} =0 .\]
\item If $f:X \rmap Y$ is continuous, $\A \in Sh(X)$ then:
\[\gc(Y;f_{\, !}\A) \simeq \gc(X;\A) \ \ , \ \  (f_{\, !}\A)_{y} \simeq \gc(f^{-1}\ps y \pd ;\A)  \ \ \ \forall\ y \in Y.\]
\end{enumerate}
}
\end{num}

\begin{num}
\label{nott}
\emph{
{\bf Notation}: Let $\A \in Sh(X) , \B \in Sh(Y)$ be c-soft sheaves. Most of the maps  we are going to deal with are of type:
\[ (\al,f)_{*}:\gc(X;\A) \rmap \gc(Y,\B) \ , \  (\ps \al,f \pd _{*}u) \ps  y \pd = \sum_{x \in f^{-1}\ps y \pd} \al_{x}(u\ps x \pd ) \in \B_{y} \ \ \ ( y \in Y) ,\] 
for some  \'etale map $f:X \rmap Y $ , and some morphism of sheaves $\al : \A \rmap f^{*}\B$. In other words, $(\al,f)_{*}$ is the composition:
\[ \gc(X;\A) \stackrel{\alpha_{*}}{\rmap} \gc(X; f^{*}\B) \stackrel{\ref{basicsgc}.3}{\simeq} \gc(Y, f_{\,!}f^{*}\B) \stackrel{int_{*}}{\rmap} \gc(Y,\B) ,\]
where $int: f_{\,!}f^{*}\B \rmap \B$ is, at the stalk at $y \in Y$ (see also \ref{basicsgc}.3):
\[ \bigoplus_{x \in f^{-1}\ps y\pd}\B_{y} \rmap \B_{y}, \ \ \ \sum_{x \in f^{-1}\ps y\pd}b_{x}x \mapsto \sum b_{x} .\]
Rather than writing the explicit formula for $\ps \al, f \pd _{*}$, we prefer to briefly indicate the maps $f$ and $\alpha$ using the notation:
\[ \gc (X,\A ) \rmap \gc (Y,\B )  \ , \  ( a | x ) \mapsto ( \al \ps a \pd | f \ps x \pd )\]
($x\in X, a \in \A _{\ x})$ }.
\end{num}

\begin{num}
\label{our}
\emph{
{\bf $\C _{c} ^{\infty}$ for non-Hausdorff manifolds}(\cite{CrMo}): If $M$ is a manifold, not necessarily Hausdorff, we put $\C _{c}^{\infty}(M):= \gc(M, \C _{M}^{\infty})$ where $\C _{M} ^{\infty}$ is the sheaf of smooth (complex valued) functions on $M$. From the Mayer-Vietoris sequence, we have an alternative description of $\C _{c}^{\infty}(M)$, as the cokernel of:
\[ \bigoplus _{U, V} \C _{c}^{\infty} (U \cap V) \rmap \bigoplus _{U} \C _{c}^{\infty} (U)  \ \ \ \ (U \in \U ) ,\]
where \U\ is a  Hausdorff open cover of $M$, and $\C _{c}^{\infty} (U)$ has the usual meaning for Hausdorff $U$ (the definition does not depend on the cover $\U$). In the same way define $\Omega_{c}^{k}(M)$ by using the sheaf $\Omega_{M}^{k}$ of (complex) $k$-forms on M. These are the natural objects with the property that constructions performed in coordinate charts patch globally. For instance, this is the case for de Rham differential: the classical one does exist locally and defines a c-soft resolution of $\mathbb{C}$:
\[ \mathbb{C} \rmap \Omega_{M}^{0} \stackrel{d_{dRh}}{\rmap}  \Omega_{M}^{1} \stackrel{d_{dRh}}{\rmap} . \ . \ .  \ \ .\]
By applying $\Gamma_{c}(M, -)$ to it, we get a globally defined de Rham differential, and a (compactly supported) de Rham complex:
\[ 0 \rmap \Omega_{c}^{0}(M) \stackrel{d_{dRh}}{\rmap}  \Omega_{c}^{1}(M) \stackrel{d_{dRh}}{\rmap} . \ . \ . \ \ ,\]
computing $H_{c}^{*}(M; \mathbb{C})$. \\
Denoting by $\C_{c,Connes}^{\infty}(M)$ Connes' definition of the vector space of ``compactly supported smooth functions'' (see section 6 in \cite{Connes}), it is not difficult to see there is a surjection:
\[ \C _{c}^{\infty}(M) \rmap \C_{c,Connes}^{\infty}(M) ,\]
which is not injective in general (of course, if $M$ is Hausdorff, the two definitions coincide and have the usual meaning).
}
\end{num}


\subsection{Homology and Cohomology of \'Etale Groupoids}

\begin{num}
\label{sheaves}
\emph{
{\bf \G -sheaves}(\cite{SGA,Mo4}): Let \G\ be an \'etale groupoid. A \G-sheaf is a sheaf \A\ on the space \nG{0}, equipped with a continuous right action of \G. This means that for any arrow $g:c \rmap d$ in \G, there is a morphism between stalks $\A_{d} \rmap \A_{c}, a \mapsto ag$, satisfying the usual identities for an action. Viewing \A\ as an \'etale space $\A \rmap \nG{0}$ (i.e. \A\ is the disjoint union of all stalks $\A_{c}$ with the germ topology) it gives a map $m:\A \times_{\nG{0}} \G \rmap \A, \ps a,g \pd \mapsto ag$; the continuity of the action means that $m$ is continuous.
\\ There is an obvious notion of morphisms of \G-sheaves; this leads to the category $Sh(\G)$ of \G-sheaves (of complex vector spaces cf. \ref{ass}). If $P:\G \rmap \Ha$ is a morphism of \'etale groupoids, it induces a functor $P^{*}: Sh(\Ha) \rmap Sh(\G)$. This construction is natural; in particular, a Morita equivalence of \'etale groupoids $\G \rmap \Ha$ induces an equivalence between their categories of sheaves; see 1.1, 2.2 in \cite{Mo1}, 4.2 in \cite{Ha3}.
}
\end{num}

\begin{num}
\label{sigma}
\emph{
{\bf Examples}: 
\begin{enumerate}
\item The constant sheaf $\mathbb{C}$, with the trivial action is a \G-sheaf. The tensor product of two \G-sheaves is naturally a \G-sheaf.
\item The \'etale map $\alpha_{n}:\nG{n}=\{ \ps g_{1}, ... , g_{n} \pd \in \G^{\ n} : s \ps g_{i} \pd = t \ps g_{i+1} \pd \} \rmap \nG{0} , \ps g_{1}, ... , g_{n} \pd \mapsto t \ps g_{1} \pd $ induces a \G-sheaf, denoted $\mathbb{C}[\nG{n}]$ (cf 2.2 in \cite{Ha3} ). Its stalk at $c \in \nG{0}$ is $\mathbb{C}[\alpha_{n}^{-1} (c)]$ and the action is given by $\ps g_{1}, ... , g_{n} \pd g=\ps g^{-1}g_{1}, g_{2}, ... , g_{n} \pd$.
\item If \G\ is a smooth \'etale groupoid than $\A=\C_{\nG{0}}^{\infty}$ is a \G-sheaf, with the action described as follows. Let $g:c \rmap d$ be an arrow in \G. Choose $V$ a neighborhood of $g$ in $\G$ such that $s|_{V} , t|_{V}$ are diffeomorphisms and put $\sigma_{g,V}:=(t|_{V}) (s|_{V})^{-1}:(s \ps V \pd , c) \rmap (t \ps V \pd , d)$. Its germ at $c$, denoted $\sigma_{g}: (\nG{0}, c) \rmap (\nG{0}, d)$  doesn't depend on the choice of $V$; it induces     multiplication by $g: \A_{d} \rmap \A_{c} , \f \mapsto \f \compose \sigma_{g}$.
\item As in the previous example, $\Omega_{\nG{0}}^{k}=$the sheaf of(complex valued) $k$-forms on $\nG{0}$ is a \G-sheaf  for all $k \geq 0$.
\end{enumerate}
}
\end{num}

\begin{num}
\label{coh}
\emph{
{\bf Cohomology}: Let \G\ be an \'etale groupoid. The category $Sh(\G)$ is an abelian category with enough injectives (\cite{SGA}); due to this fact, the cohomology of sheaves on \G\ can be defined and used via elementary homological algebra tools; using standard resolutions, it can be computed by some kind of standard bar-complexes; see e.g. \cite{Ha3,Mo4,SGA} (a brief review will be given in \ref{cohomol}).
}
\end{num}
\ \ \ More generally, one can define the bi-functors $Ext_{\G}^{*}(-,-): Sh(\G) \times Sh(\G) \rmap \underline{Vs}$ \  ($\underline{Vs}$ denotes the category of complex vector spaces) by:
\[ Ext_{\G}^{*}(\A ,\B )=R^{*}Hom_{Sh\ps \G \pd}( \A, - )( \B ) ,\]
with the particular case $Ext_{\G}^{*}(\mathbb{C},-)=H^{*}(\G ; -)$ ($R^{*}$ stands for the right derived functors \cite{Hi,We1}).
\par Homological algebra provides us an alternative description of the vector spaces $Ext_{\G}^{p}(\A ,\B )$ by means of Yoneda extensions (see \cite{MaL}, \cite{We1}). For $p \geq 1$, the elements in $Ext_{\G}^{p}(\A ,\B )$ are represented by p-extensions of \A\ by \B\ i.e. exact sequences in $Sh(\G)$: 
\[ u: \ \ \ 0 \rmap \B \rmap X_{1} \rmap .\ .\ . \rmap X_{p} \rmap \A \rmap 0 .\]
The equivalence relation is generated by: $u \simeq \tilde{u}$ whenever there exists a morphism of complexes $u \mapsto \tilde{u}$. According to this, there is a simple description of the cup-product  $Ext_{\G}^{p}(\A ,\B ) \times Ext_{\G}^{q}(\B ,\C ) \rmap  Ext_{\G}^{p+q}(\A ,\C )$ for $\A, \B, \C \in Sh(\G)$ as the concatenation of exact sequences (for the general setting of abelian categories see \cite{MaL}, pp. 82-87, \cite{We1}, pp. 76-80).
\par We have also a simple description of the cap-products: for any homological $\delta$-functor (2.1.1 in \cite{We1}) $L_{*}:Sh(\G) \rmap \underline{Vs}$  \ there are cap-product maps $ \xymatrix{
L_{n}(\A) \times Ext_{\G}^{p}(\A ,\B ) \ar[r]^-{\cap} & L_{n-p}(\B)}$ ($p \geq 0 , n , p$ integers, $\A,\B \in Sh(\G)$). If $p=0$, $Ext_{\G}^{0}(\A ,\B )=Hom_{Sh\ps \G \pd}( \A, - )( \B )$ and $\cap$ is the covariance of $L_{n}$. If $p=1 ,\ u \in Ext_{\G}^{1}(\A ,\B )$ the cap-product by $u$ : $-\cap u : L_{n}(\A) \rmap L_{n-1}(\B)$ is the boundary of the long exact sequence associated to (any) short exact sequence $0 \rmap \B \rmap X \rmap \A \rmap 0$ representing $u$. If $p \geq 2$, we iterate the case $p=1$.
\par For any $\G$-sheaf $\A$ there is an obvious morphism $Ext_{\G}^{p}(\mathbb{C},\mathbb{C} ) \rmap Ext_{\G}^{p}(\A ,\A )$ (tensoring by \A); we get in particular an action of the cohomology on any homological $\delta$-functor: \[ L_{n}(\A) \times H^{p}(\G,\mathbb{C}) \rmap L_{n-p}(\A) \ \ , \ \ \A \in Sh(\G) \ \ .\]

\begin{num}
\label{classif}
\emph{
{\bf Connection with the classifying space}: Any \G-sheaf \es \ gives rise to a sheaf $\tilde{\es}$ on the classifying space $B\G$ and there are isomorphisms $H^{*}(\G;\es) \simeq H^{*}(B\G;\tilde{\es})$. This was conjectured by Haefliger and proved in \cite{Mo2}.
}
\end{num}

\begin{num}
\label{barcomp}
\emph{
{\bf Bar-complexes}: Let \G\ be an \'etale groupoid. Denote:
\[ \nG{n}= \{ \ps g_{1}, ... , g_{n} \pd \in \G^{\ n} : \ s \ps g_{i} \pd = t \ps g_{i+1} \pd \ , \ \forall \ 1 \leq i \leq n-1 \} .\]
\ \ \ For any $\A \in Sh(\nG{0})$, we consider for each n the pull-back of $\A$ to \nG{n} along $\nG{n} \rmap \nG{0} , \ps g_{1}, ... , g_{n} \pd \mapsto t \ps g_{1} \pd$. These sheaves on the spaces $\nG{n}$ are still denoted by \A\ (it will be clear from the context on which space they are considered). For any c-soft sheaf $\A \in Sh(\G)$ (c-softness of $\G$-sheaves will always mean c-softness on $\nG{0}$) , define the bar-complex $B_{\dt}(\G;\A)$ as the chain complex associated to the simplicial vector space $n \mapsto \gc (\nG{n};\A)$ with the structure maps:  \[ d_{i} \ps a \ | \ g_{1},. . ., g_{n} \pd = \left \{ \begin{array}{ll} 
                                              \ps ag_{1} \ | \ g_{2}, . . . , g_{n} \pd & \mbox{if $i=0$} \\
                                              \ps a \ | \ g_{1}, . . . , g_{i}g_{i+1}, . . . , g_{n} \pd & \mbox{if $1 \leq i \leq n-1$} \ \ ,\\
                                              \ps a \ | \ g_{1}, ... , g_{n-1} \pd & \mbox{if $i=n$}
                                                        \end{array}
                                            \right.  \]
\[ s_{i} \ps a \ |  \ g_{1}, ... , g_{n} \pd = \ps a\  |\  . . ., g_{i}, 1, g_{i+1}, . . . \pd ,\]
(see  \ref{nott} for notations). 
}
\end{num}

\begin{num}
\label{homol}
\emph{
{\bf Homology}(\cite{CrMo}): Despite to the fact that $Sh(\G)$  does not have enough projectives in general, there is quite an  obvious way to define the homology of sheaves on \'etale groupoids by using bar-complexes and c-soft resolutions. Looking at the classical homology of groups, one can state and prove "usual results" like long exact sequences, spectral sequences, etc. We review some definitions and properties we need, referring for a more detailed description to \cite{CrMo}. The analogous properties for cohomology are well known (\cite{Ha3,SGA,Mo4}).
}
\end{num}
\par If $\A \in Sh(\G)$ is c-soft, define $H_{*}(\G,\A)$ as the homology of $B_{\dt}(\G,\A)$. If $\A_{\dt}$ is a chain complex of c-soft \G-sheaves, define $\mathbb{H}_{*}(\G,\A_{\dt})$ as the homology of the double complex  $B_{\dt}(\G,\A_{\dt})$. If $\es \in Sh(\G)$, take $\es \rmap \A^{0} \rmap \A^{1} \rmap .\ .\ .\ $ a c-soft resolution in $Sh(\G)$ and define $H_{*}(\G,\es)$ as $\mathbb{H}_{*}(\G,\tilde{\A}_{\dt})$ where $\tilde{\A}_{k}=\A^{-k} \ \ \forall \ k$-integer. The notion of Cartan-Eilenberg resolutions (5.7.9 in \cite{We1}) by c-soft sheaves carries over to $Sh(\G)$ and one can define $\mathbb{H}_{*}(\G,\es_{\dt})$ for any complex of \G-sheaves $\es_{\dt}$. The existence of resolutions and that the result does not depend on the choices we make are proved in \cite{CrMo}. 

\begin{num}
\label{uit}
\emph{ From the previous definition and \ref{basicsgc}.1 we see that any morphism $\alpha: \es_{\dt} \rmap \es^{\, '}_{\dt}$ of complexes of $\G$-sheaves which is a quasi-isomorphism (of complexes of sheaves on $\nG{0}$), induces an isomorphism $\mathbb{H}_{*}(\G; \es_{\dt}) \rmap \mathbb{H}_{*}(\G; \es_{\dt}^{\, '})$.
}
\end{num}

\begin{num}
\label{lili}
\emph{
{\bf Example}: If $\G=G$ is a group, we get the usual homology of groups (\cite{Br}). If $\G=X$ is a space we get $H_{k}(\G;-)=H_{c}^{-k}(X;-)$. In general, $H_{*}(\G;-)$ lives in degrees $* \geq -cohdim(\nG{0})$.
}
\end{num}

\begin{num}
\label{leseq}
\emph{
{\bf The long exact sequence}(\cite{CrMo}): For any short exact sequence $0 \rmap \es^{'} \rmap \es \rmap \es^{''} \rmap 0$ in $Sh(\G)$ there is a long exact sequence of vector spaces: 
\[ .\ .\ .\ \rmap H_{n}(\G;\es^{'}) \rmap H_{n}(\G;\es) \rmap H_{n}(\G;\es^{''}) \rmap H_{n-1}(\G;\es^{''}) \rmap .\ .\ . \]
}
\end{num}
and it is natural with respect to morphisms in $Sh(\G)$. In other words, $H_{*}(\G;-)$ is a homological $\delta$ -functor on $Sh(\G)$.  
\par In particular we get (cf. \ref{coh}) the cap-products $H_{n}(\G;\A) \times Ext^{p}_{\G}(\A;\B) \rmap H_{n-p}(\A;\B)$ and $H_{n}(\G;\A) \times H^{p}(\G;\mathbb{C}) \rmap H_{n-p}(\G;\A)$.

\begin{num}
\label{spseq}
\emph{
{\bf Basic spectral sequences}(\cite{CrMo}): For any bounded below chain complex $\es_{\dt}$ in $Sh(\G)$ there are two spectral sequences:
\[ E_{p,q}^{2}=H_{p}(H_{q}(\G;\es _{\dt})) \Longrightarrow \mathbb{H}_{p+q}(\G;\es_{\dt}) ,\]
\[ E_{p,q}^{2}=H_{p}(\G;\tilde{H}_{q}(\es_{\dt})) \Longrightarrow \mathbb{H}_{p+q}(\G;\es_{\dt}) ,\]
(compare to 5.6 in \cite{Br}). In the second spectral sequence, $\tilde{H}_{q}(\es_{\dt}) \in Sh(\G)$ denotes the homology sheaf of $\es_{\dt}$. Given $\es_{\dt}$ we define an element:
\[ u_{q}(\es_{\dt}) \in Ext_{\G}^{2}(\tilde{H}_{q}(\es_{\dt}),\tilde{H}_{q+1}(\es_{\dt}) )\]
as follows. Denote $\Ka_{q} = Ker(\es_{q} \rmap \es_{q-1})$. The exact sequence $0 \rmap \Ka_{q+1} \rmap \es_{q+1} \rmap \Ka_{q} \rmap \tilde{H}_{q}(\es_{\dt}) \rmap 0$ defines an element $\epsilon_{q} \in Ext_{\G}^{2}(\tilde{H}_{q}(\es_{\dt}); \Ka_{q+1} )$. Using the projection $\pi : \Ka_{q+1} \rmap \tilde{H}_{q+1}(\es_{\dt})$ and the covariance of $Ext$ in  the second variable, put $u_{q}(\es_{\dt})=\pi_{*}(\epsilon_{q})$.
}
\end{num}

\begin{num}
\label{bdr}
{\bf Lemma}: The $d^{2}-$ boundaries of the second spectral sequence are:
\[ d_{p,q}^{2}= - \cap u_{q}(\es_{\dt}) : H_{p}(\G;\tilde{H}_{q}(\es_{\dt})) \rmap  H_{p-2}(\G;\tilde{H}_{q+1}(\es_{\dt})) .\]
\end{num}

\emph{proof}: The spectral sequence is defined as follows (\cite{CrMo}). Take $\mathbb{C} \rmap \ii^{\ \dt}$ to be a bounded c-soft resolution of $\mathbb{C}$ in $Sh(\G)$ and put $\ii_{k}=\ii^{-k}$. Then $\ii^{\ \dt} \otimes \es_{\dt}$ is a Cartan-Eilenberg resolution of $\es_{\dt}$ so $\mathbb{H}_{*}(\G; \es_{\dt})$ is computed by the triple complex $B_{\dt}(\G; \ii_{\ \dt} \otimes \es_{\dt})$, or, equivalently, by the double complex $C_{\dt \dt}$ defined by:
\[ C_{p, q} = \bigoplus_{p_{1}+p_{2}=p} B_{p_{1}}(\G; \ii_{p_{2}} \otimes \es_{q}) .\]
One of its spectral sequences has $E_{p, q}^{1}=H_{q}(C_{p, \dt})= \bigoplus_{p_{1}+p_{2}=p} B_{p_{1}}(\G; \ii_{p_{2}} \otimes \tilde{H}_{q}(\es_{\dt}))$ and $E_{p, q}^{2}=H_{p}(\G; \tilde{H}_{q}(\es_{\dt}))$. From the general description of the boundaries of the spectral sequence associated to a double complex, our boundary is described by a "zig-zag" of length two (see \cite{BoTu}, page 164). This corresponds to two boundaries of long exact sequences i.e. to a cap-product by a $Ext^{2}$-class; the relation $d_{p,q}^{2}= - \cap u_{q}(\es_{\dt})$ becomes a standard checking inside the spectral sequence.

\begin{num}
\label{HS}
\emph{
{\bf Hochschild-Serre spectral sequence}(\cite{CrMo}): If $\f: \G \rmap \Ha$ is a continuous morphism between  \'etale groupoids, $\es \in Sh(\G)$, there is a spectral sequence:
\[ E_{p,q}^{2}=H_{p}(\Ha; L_{q}\f_{\,!} (\es)) \Longrightarrow H_{p+q}(\G; \es) .\]
Here $L_{q}\f_{\,!} (\es) \in Sh(\G)$ is a sheaf with stalks:
\[ L_{q}\f_{\,!} (\es)_{d}=H_{q}(d/\f; \omega_{d}^{*} \es) \ \ \ , \forall \ d \in \nH{0} \]
(compare to the similar result for cohomology, \cite{Mo4}, pp 15-16).
}
\end{num}

\ \ \ For $q=0$, denote $L_{0}\f_{\,!}=\f_{\,!}$. As in the case of spaces (\cite{We1,Iv}), the spectral sequence is a consequence of the second spectral sequence in \ref{spseq} and of an equality: $H_{*}(\G; \es)= \mathbb{H}_{*}(\Ha ; \el \f_{\,!}(\es))$. Here $\el \f_{\,!}(\es)$ is a chain complex in $Sh(\Ha)$ with $\tilde{H}_{q}(\el \f_{\,!}(\es))=L_{q}\f_{\,!}(\es)$. In particular, the second boundary is of type:
\[ d_{p,q}^{2}= - \cap u_{q} : E_{p,q}^{2} \rmap E_{p-2,q+1}^{2} \ \ \ \mbox{with} \ u_{q} \in Ext_{\Ha}^{2}(L_{q}\f_{\, !}\es;L_{q+1}\f_{\, !}\es) .\]

\begin{num}
\label{mor}
\emph{
{\bf Morita invariance}(\cite{CrMo}): Any morphism $P: \G \rmap \Ha$ which is \'etale (in the sense that $s_{P}$ is \'etale), induces a morphism in homology $P_{*}: H_{*}(\G;P^{*}\B) \rmap H_{*}(\Ha;\B)$ for any $\B \in Sh(\G)$. This construction is natural with respect to \'etale morphisms; in particular, $H_{*}(\G;-)$ is Morita invariant and any essential equivalence $\f :\G \rmap \Ha$ (see \ref{eseq}.1) induces an isomorphism $H_{*}(\G; \f^{*}\B) \simeq H_{*}(\Ha; \B)$.
}
\end{num}

\begin{num}
\label{poincare}
\emph{
{\bf Duality}: For any \'etale groupoid \G\ with \nG{0} a topological manifold of dimension $n$, there is a Poincar\'e-duality isomorphism $H^{k}(\G;or) \simeq H_{k-n}(\G;\mathbb{C})^{\vee} \ \ ,\forall \ k \geq 0$. Here $or \in Sh(\G)$ is the complex orientation sheaf of \nG{0} (\cite{Iv})  with the obvious \G-action, and $\vee$ denotes the vector-space dual. When $\G=M$ is a topological manifold, this is the usual Poincar\'e duality. See \cite{CrMo} for the full Verdier-duality for \'etale groupoids.
}
\end{num}

\begin{num}
\label{Etcat}
\emph{
{\bf \'Etale categories}: Without any changes, the definitions and the basic properties we have described so far in this sub-section work equally well for \'etale categories (i.e. small categories \G, with topologies on $\nG{0}, \nG{1}$ such that all the structure maps are \'etale) with only one exception: Morita invariance (to prove Morita invariance for \emph{\'etale categories} we have to take care of what a morphism between \'etale categories is; see 5.4 in \cite{Mo1} or \cite{Mo4}).
Because of this we need another technical tool when we deal with \'etale categories (for a more general result, see \cite{CrMo}):
}
\end{num}

\begin{num}
\label{inv}
{\bf Lemma and definition}: Let $\G$ and $\Ha$ be \'etale categories. A continuous functor $\f: \G \rmap \Ha$ is called a strong deformation retract of $\Ha$ if there is a continuous functor $\psi: \Ha \rmap \G$ (called retraction) and a continuous natural transformation of functors $F: \f \compose \psi \rmap Id_{\Ha}$ (called strong deformation retraction) such that $\psi \compose \f = Id_{\G}, F(\f\ps c \pd)=id_{\f \ps c \pd}$ for all $c \in \nG{0}$ and the maps $\f, \psi, F$ are \'etale.
\par In this case, for any $\Ha$-sheaf $\A$, $\f$ induces an isomorphism:
\[  H_{*}(\G ; \f^{*}\A) \tilde{\rmap} H_{*}(\Ha; \A) .\]
\end{num}
\emph{proof}: Denote $\f \compose \psi$ by $l$. Let $\Phi$ the map induced by $\f$:
\[ \Phi\ : B_{\dt}(\G; \f^{*}\A) \rmap B_{\dt}(\Ha; \A), \ \ \Phi (a| g_{1}, . . . , g_{n})= (a| \f \ps g_{1}\pd , . . . , \ps g_{n} \pd) .\]
Since $\f$ is \'etale, $\f^{*}: Sh(\nH{0}) \rmap Sh(\nG{0})$ preserves c-softness (\cite{Iv}) so it is enough to prove that $\Phi$ is a homotopy equivalence of chain complexes, when $\A$ is c-soft. Define a chain map:
\[ \Psi : B_{\dt}(\Ha; \f^{*}\A) \rmap B_{\dt}(\G; \A), \ \ \Psi (a| h_{1}, . . . , h_{n})= (aF(t\ps h_{1} \pd)| \psi\ps h_{1}\pd , . . . ,\psi\ps h_{n}\pd) .\] 
We have $\Psi\ \compose \Phi\ =Id$ and $\Phi\ \compose \Psi$ is homotopic to $Id$ by the following homotopy:
\[ h: B_{\dt}(\Ha; \f^{*}\A) \rmap B_{\dt\ +1}(\Ha; \f^{*}\A), \ \ \ h= \sum_{i=0}^{n} (-1)^{i}h_{i} ,\]
\[ h_{i}(a | h_{1}, . . . , h_{n})=(a | h_{1}, . . . , h_{i}, F(s\ps h_{i}\pd) , l\ps h_{i+1}\pd, . . . , l\ps h_{n}\pd).  \]


\section{Cyclic Homologies of Sheaves on \'Etale Groupoids}


\subsection{Cyclic Objects}
\ \ \ Recall (\cite{Co1,FeTs,Lo}) some basic definitions concerning cyclic objects in an abelian category \M.

\begin{num}
\label{mixed}
\emph{
{\bf Mixed Complexes}: A mixed complex in \M\ is a family $\{X_{n}: n \geq 0 \}$ of objects in \M, equipped with maps of degree $-1$, $b: X_{n} \rmap X_{n-1}$, maps of degree $1$, $B: X_{n} \rmap X_{n+1}$, satisfying $b^{2}=B^{2}=bB+Bb=0$. A mixed complex gives rise to a first quadrant double complex in \M:
}
\end{num}
\[  \xymatrix{
\ar[d]_{b}                     &   \ar[d]_{b}                      &   \ar[d]_{b}                         & \ar[d]_{b} \\
X_{3} \ar[d]_{b}               &   X_{2} \ar[l]_{B} \ar[d]_{b}     &   X_{1} \ar[l]_{B} \ar[d]_{b}        & X_{0} \ar[l]_{B}  \\
X_{2} \ar[d]_{b}               &   X_{1} \ar[d]_{b} \ar[l]_{B}     &   X_{0} \ar[l]_{B} \\
X_{1} \ar[d]_{b}               &   X_{0} \ar[l]_{B} \\
X_{0}
}
\]
\ \ \ It (or its total complex) is denoted $(X_{\dt},B,b)$. The Hochschild and cyclic homology of $X_{\dt}$ are defined by $HH_{*}(X_{\dt}):=H_{*}(\ps X_{\dt},b \pd), \ HC_{*}(X_{\dt}):=H_{*}(\ps X_{\dt},B,b\pd)$. There is a short exact sequence :
\[ 0 \rmap (X_{\dt},b) \stackrel{I}{\rmap} (X_{\dt},B,b) \stackrel{S}{\rmap} (X_{\dt},B,b)[-2] \rmap 0 ,\]
where $I$ is the inclusion on the first column and $S$ is the shifting. Standard homological algebra implies there is a long exact sequence :
\[ . \ . \ . \stackrel{B}{\rmap} HH_{n}(X_{\dt}) \stackrel{I}{\rmap} HC_{n}(X_{\dt}) \stackrel{S}{\rmap} HC_{n-2}(X_{\dt}) \stackrel{B}{\rmap} HH_{n-1}(X_{\dt}) \stackrel{I}{\rmap} . \ . \ . \]
(the SBI-sequence of $X_{\dt}$). Using the shift operator, the periodic cyclic homology of $X_{\dt}$ is defined by $HP_{*}(X_{\dt}):=H_{*}(\lim_{r}(X_{\dt},B,b)[-2r])$.

\begin{num}
\label{ciclic}
\emph{
{\bf Cyclic Objects}: Usually, mixed complexes are made out of cyclic objects. Let $1 \leq r \leq \infty$. An $r$-cyclic object in \M\ is a contravariant functor $X:\la_{r} \rmap \M$ from the generalized Connes  category $\la_{r}$ (A2 in \cite{FeTs}). That means a simplicial object $\{ X_{\dt},d_{\dt},s_{\dt}\}$ in \M\ together with morphisms $t_{n}:X_{n} \rmap X_{n}$ such that:
\[ d_{i}t_{n}= \left\{ \begin{array}{ll}
                                      t_{n-1}d_{i-1} & \mbox{if $i \neq 0$} \\ 
                                      d_{n} & \mbox{if $i=0 $}
                        \end{array}
                \right.
\ , \ \   s_{i}t_{n}=\left\{ \begin{array}{ll}
                                      t_{n+1}s_{i-1} & \mbox{if $i \neq 0$} \\
                                      t_{n+1}^{2}s_{n} & \mbox{if $i=0$}
                        \end{array}
                \right.  ,  \]
}
\end{num}
and the cyclic relation $t_{n}^{r\ps n+1 \pd}=1$ holds in the case $r \neq \infty$.
\par We define $b^{'}, b:X_{n} \rmap X_{n-1}, b^{'}=\sum_{j=0}^{n-1}\ps -1 \pd^{j}d_{j}, \ b=b^{'}+ \ps -1 \pd^{n}d_{n}, \ s_{-1}=s_{n}t_{n+1}$ the extra degeneracy (which gives a contraction for $b^{'}$), $\tau_{n}=\ps -1 \pd^{n}t_{n+1}$. If $r \neq \infty$, let $N=\sum_{j=0}^{\ps n+1 \pd r -1}\tau_{n}^{j}, \ B=\ps 1-\tau_{n} \pd s_{-1}N$. Then $(X_{\dt},b,B)$ is a mixed complex; its homologies are denoted $HH_{*}(X_{\dt}), HC_{*}(X_{\dt})$. They are in fact the homologies of $X_{\dt}$ as a contravariant functor $\la_{\infty}\rmap \M$ respectively $\la_{r} \rmap \M$ (cf. A3.2 in \cite{FeTs}).
\par For $r=\infty$ we put $HH_{*}(X_{\dt}):=H_{*}((X_{\dt},b))$.
\par For $r=1$, the category $\la_{1}$ is denoted $\la$ and the $1$-cyclic objects are called cyclic objects.

\begin{num}
\label{alfaex}
\emph{
{\bf Examples}: The basic example (\cite{FeTs}, \cite{Ni1}) is the $\infty$-cyclic vector space $A_{\al}^{\natural}$ associated to an unital algebra $A$, endowed with an endomorphism $\al: A_{\al}^{\natural}\ps n \pd:=A^{\otimes \ps n+1 \pd}$,
\[ d_{i}( a_{0}, . \ . \ . , a_{n})=\left\{ \begin{array}{ll}
                                                       ( a_{0}, . \ . \ .\ ,a_{i}a_{i+1}, . \ . \ .\ ,a_{n} )& \mbox{if $0 \leq i \leq n-1$} \\
                                                       ( \al \ps a_{n} \pd a_{0}, a_{1}, . \ . \ . \ , a_{n-1} )& \mbox{if $i=n$}
                                                \end{array}
                                        \right. ,\]
\[ t ( a_{0}, . \ . \ . \ , a_{n}) = ( \al \ps a_{n} \pd , a_{0}, a_{1}, . \ . \ .\ , a_{n-1} ) \ ,\]  
\[ s_{i} ( a_{0}, . \ . \ . \ , a_{n} ) = ( \ . \ . \ . \ , a_{i}, 1, a_{i+1}, . \ . \ .\ ) .\] 
}
\end{num}
\ \ \ Its Hochschild homology is denoted $HH_{*}(A,\al)$. If $\al$ is of order $r \neq \infty$ than $A_{\al}^{\natural}$ is an $r$-cyclic vector space; denote by $HC_{*}(A;\al), HP_{*}(A;\al)$  the corresponding homologies. 
\par For $\al=id$, we simplify the notations to $A^{\natural}, HH_{*}(A), HC_{*}(A) , HP_{*}(A)$.These can be defined more generally, for non-unital algebras $A$, by using the $\ps b,b^{'} \pd $-complex (see \cite{Co3,Lo}):
\[ \xymatrix{
 \ar[d]_{b} &  \ar[d]_{b^{'}} &  \ar[d]_{b} &  \ar[d]_{b^{'}} \\
A^{\otimes 3} \ar[d]_{b} & A^{\otimes 3} \ar[d]_{b^{'}} \ar[l]_{1-\tau } & A^{\otimes 3} \ar[d]_{b} \ar[l]_{N} & A^{\otimes 3} \ar[d]_{b^{'}} \ar[l]_{1-\tau } & . \ . \ . \ \ar[l]_{N} \\
A^{\otimes 2} \ar[d]_{b} & A^{\otimes 2} \ar[d]_{b^{'}} \ar[l]_{1-\tau } & A^{\otimes 2} \ar[d]_{b} \ar[l]_{N} & A^{\otimes 2} \ar[d]_{b^{'}} \ar[l]_{1-\tau } & . \ . \ . \ \ar[l]_{N} \\
A & A \ar[l]_{1-\tau } & A \ar[l]_{N} & A \ar[l]_{1-\tau } & . \ . \ . \ \ar[l]_{N}
 } \]
\ \ \ When $A$ is a locally convex algebra, we use the projective tensor product (\cite{Gr}).

\begin{num}
\label{topology}
\emph{
{\bf Example}: If $\A=\C_{M}^{\infty}$ is the sheaf of smooth functions on a manifold $M$, we define $\A^{\natural}$ as in \ref{alfaex} by taking into account the topology, i.e. $\A^{\natural}\ps n \pd :=\Delta_{n+1}^{*}( \A ^{\boxtimes \ps n+1 \pd } )$ where $\A^{\boxtimes n}:=\C_{M^{n}}^{\infty} \in Sh(M^{n})$ and $\Delta_{n}: M \rmap M^{n}$ is the diagonal map. Keeping the same formulas as in \ref{alfaex} with $\al=id$, $\A^{\natural}$ is a cyclic object in $Sh(M)$.
}
\end{num}
\ \ \ If $\f : (M,x_{0}) \rmap (M,x_{0})$ is the germ of a smooth function, it induces an algebra endomorphism $\A_{x_{0}} \rmap \A_{x_{0}}$. Using the  formulas from \ref{alfaex} we get an $\infty$-cyclic vector space. Denote its $b$-boundary by $b_{\f}$. The following lemma belongs to folklore (a proof is included in the proof of 5.1 and 5.2 in \cite{BrNi}):

\begin{num}
\label{stability}
{\bf Lemma}: We call $\f : (M,x_{0}) \rmap (M,x_{0})$ stable if it satisfies the following three conditions:
\par (i)  $M^{\f}:=\{ x \in dom(\f ): \f \ps x \pd =x \}$ is locally around $x_{0}$ a submanifold of $M$;
\par (ii) $T_{x_{0}}M^{\f}=Ker( \ps d \f \pd _{x_{0}} - Id)$;
\par (iii)$T_{x}M$ splits into a direct sum of $T_{x_{0}}M^{\f}$ and a \f-invariant subspace.
\par Any germ \f\ which preserves a metric around $x_{0}$ (in particular any germ of finite order) is stable. If \f\ is stable, then:
\[ ( (\C_{M}^{\infty})_{x_{0}}^{\natural} , \ b_{\f} \ ) \rmap \ ( (\C_{M^{\f}}^{\infty})_{x_{0}}^{\natural}, b_{id} \ )  \rmap \ ( (\Omega_{M^{\f}}^{*})_{x_{0}} , \ 0 \ ) , \]
\[ \ \ps a_{0}, .  .  . , a_{n} \pd  \mapsto \ps a_{0}|_{M^{\f}}, .  .  . , a_{n}|_{M^{\f}} \pd \mapsto a_{0}\ da_{1} \ . \ . \ . \ da_{n} |_{M^{\f}} \]
are quasi-isomorphisms of complexes of vector spaces.
\end{num}


\subsection{Cyclic \G\ -sheaves}

\begin{num}
\label{ciclicsheaves}
{\bf Definition}: Let \G\ be an \'etale groupoid. By a cyclic \G-sheaf we mean a cyclic object $\A_{\dt}$ in $Sh(\G)$. Its homology \G-sheaves, as defined in \ref{ciclic}, are denoted $\widetilde{HH}_{*}(\A_{\dt}), \widetilde{HC}_{*}(\A_{\dt}), \widetilde{HP}_{*}(\A_{\dt})$.
\end{num}

\begin{num}
\label{standard}
\emph{
{\bf Examples}:
\begin{enumerate}
\item As an extension of \ref{alfaex}, any sheaf of complex algebras $\A \in Sh(\G)$ defines a cyclic \G-sheaf $\A^{\natural}$ having the stalk at $c \in \nG{0}$: $(\A^{\natural})_{c}=(\A_{c})^{\natural}$. If $\G$ is smooth and $\A=\C^{\infty}_{\nG{0}}$, then we reserve the notation $\A^{\natural}$ for the cyclic \G-sheaf defined by taking into account the topology (see \ref{topology}, \ref{sigma}).
\item The standard resolution of $\mathbb{C}$ on $Sh(\G)$ (\cite{CrMo,Ha3,Mo4}):
\[ . \ . \ . \ \rmap \mathbb{C}[\nG{3}] \rmap \mathbb{C}[\nG{2}] \rmap \mathbb{C}[\nG{1}] \ ( \rmap \mathbb{C} \rmap 0) ,\]
comes from the simplicial structure on the nerve of \G:
\[ d_{i}\ps  g_{0},. . ., g_{n} \pd=\left\{ \begin{array}{ll}
                                                 \ps  g_{0},. . .,g_{i}g_{i+1},. . . , g_{n} \pd & \mbox{if $0 \leq i \leq n-1$} \\
                                                 \ps  g_{0},. . ., g_{n-1} \pd & \mbox{if $i=n$}
                                            \end{array}
                                    \right. ,\]
and $s_{i}$'s inserts units (as in \ref{barcomp}). It inherits a structure of cyclic \G-sheaf by (compare to 7.4.5. in \cite{Lo}):
\[ t\ps  g_{0},. . ., g_{n} \pd=\ps g_{0}g_{1}...g_{n}, \ps g_{1}...g_{n} \pd ^{-1}, g_{1}, ... , g_{n-1} \pd .\]
\end{enumerate}
}
\end{num}

\begin{num}
\label{nustiu}
{\bf Definition}: If $\A_{\dt}$ is a cyclic \G-sheaf, define its Hochschild and cyclic hyperhomology by $HH_{*}(\G;\A_{\dt})=\mathbb{H}_{*}(\G;( \A_{\dt},b ) ), \ HC_{*}(\G;\A_{\dt})=\mathbb{H}_{*}(\G;( \A_{\dt},B,b )) \  $(compare to \cite{We2}). If $\A_{n}$ is c-soft for all $\ n$, define $HP_{*}(\G;\A_{\dt})=\mathbb{H}_{*}(\G;\lim_{r} ( \A_{\dt},B,b ) [-2r] )$ (in the general case we can define $HP_{*}$ by the same methods as in \cite{We2}).
\end{num}
This is an extension of the definition given by Loday for groups (section II in \cite{Loo}).

\begin{num}
\label{SBI}
\emph{
{\bf SBI-sequences}: From the general considerations in \ref{mixed}, there is a long exact sequence in $Sh(\G)$:
\[ . \ . \ . \stackrel{B}{\rmap} \widetilde{HH}_{n}(\A_{\dt}) \stackrel{I}{\rmap} \widetilde{HC}_{n}(\A_{\dt}) \stackrel{S}{\rmap} \widetilde{HC}_{n-2}(\A_{\dt}) \stackrel{B}{\rmap} \widetilde{HH}_{n-1}(\A_{\dt}) \stackrel{I}{\rmap} . \ . \ .  \ \ ,\]
and, using \ref{leseq}, a long exact sequence of vector spaces:
\[ . \ . \ . \stackrel{B}{\rmap} HH_{n}(\G;\A_{\dt}) \stackrel{I}{\rmap} HC_{n}(\G;\A_{\dt}) \stackrel{S}{\rmap} HC_{n-2}(\G;\A_{\dt}) \stackrel{B}{\rmap} HH_{n-1}(\G;\A_{\dt}) \stackrel{I}{\rmap} . \ . \ . \ \ .\]
}
\end{num}

\begin{num}
\label{1st}
\emph{
{\bf First spectral sequences}: Using the second spectral sequence in \ref{spseq} we get two spectral sequences with $E^{2}$-terms: $H_{p}(\G;\widetilde{HH}_{q}(\A_{\dt})) \Longrightarrow HH_{p+q}(\G;\A_{\dt})$ and $H_{p}(\G;\widetilde{HC}_{q}(\A_{\dt})) \Longrightarrow HC_{p+q}(\G;\A_{\dt})$.
}
\end{num}


\begin{num}
\label{2nd}
\emph{
{\bf Second spectral sequences}: The spectral sequences of the double complex $B_{\dt}(\G;T_{\dt})$ (see \ref{barcomp}) where $T_{\dt}=(\A_{\dt},b)$  or $(\A_{\dt},B,b)$ give two spectral sequences with $E^{2}$-terms:
\[ HH_{p}(H_{q}(\G;\A_{\dt})) \Longrightarrow HH_{p+q}(\G;\A_{\dt})  \ \  \mbox{and} \ \ HC_{p}(H_{q}(\G;\A_{\dt})) \Longrightarrow HC_{p+q}(\G;\A_{\dt}) .\]
}
\end{num} 

\begin{num}
\label{SBIARG}
{\bf Lemma}: If the morphism $f:\A_{\dt} \rmap \tilde{A}_{\dt}$ of cyclic \G-sheaves induces a quasi-isomorphism $f:(\A_{\dt},b) \rmap (\tilde{A}_{\dt},\tilde{b})$ of complexes of sheaves, then it induces isomorphisms:
\[ HH_{*}(\G;\A_{\dt}) \simeq HH_{*}(\G;\tilde{A}_{\dt}) \ , \ HC_{*}(\G;\A_{\dt}) \simeq HC_{*}(\G;\tilde{A}_{\dt}) .\]
\end{num}

\emph{proof}: this is a consequence of \ref{uit}, \ref{SBI}, \ref{1st} and comparison-theorem for spectral sequences (compare to 2.5.2 in \cite{We2}). 

\begin{num}
\label{atunits}
\emph{
Assume that $\G$ is an \'etale groupoid and $\A_{\dt}$ is a cyclic \G-sheaf such that any $\A_{n}$ is c-soft. From the definition  we see that $HH_{*}(\G;\A_{\dt})$ is computed by the bi-simplicial vector space $B_{\dt}(\G;\A_{\dt})$; so, from the Eilenberg-Zilber theorem (\cite{We1}), it is computed by its diagonal, i.e. by the simplicial vector space:
\[ \xymatrix {
\C (\G;\A_{\dt}): \  .\ .\ .\ \ar@<2ex>[r]\ar@<1ex>[r]\ar[r]\ar@<-1ex>[r] &\gc (\nG{2};\A_{2}) \ar@<1ex>[r]\ar[r]\ar@<-1ex>[r] & \gc(\nG{1};\A_{1}) \ar[r]\ar@<-1ex>[r] & \gc (\nG{0};\A_{0}) } ,\]
}
\end{num}
\[ d_{i} \ps a \ | \ g_{1},. . ., g_{n} \pd = \left \{ \begin{array}{ll} 
                                              \ps d_{0}\ps a\pd g_{1} \ | \ g_{2}, . . . , g_{n} \pd & \mbox{if $i=0$} \\
                                              \ps d_{i}\ps a\pd \ | \ g_{1}, . . . , g_{i}g_{i+1}, . . . , g_{n} \pd & \mbox{if $1 \leq i \leq n-1$} \\
                                              \ps d_{n}\ps a\pd \ | \ g_{1}, ... , g_{n-1} \pd & \mbox{if $i=n$}
                                                        \end{array}
                                            \right. ,\]  
\[ s_{i} \ps a \ |  \ g_{1}, ... , g_{n} \pd = \ps s_{i}\ps a\pd \  |\  . . ., g_{i}, 1, g_{i+1}, . . . \pd \ \ .\]
Combining the cyclic structure of $\A_{\dt}$ with the one on the nerve of \G\ (see \cite{Bu} , pag. 358), we define the following cyclic structure on $\C (\G;\A_{\dt})$:
\[ t \ps a \ | \ g_{1},. . ., g_{n} \pd = ( t\ps a \pd g_{1} ... g_{n} | \ps g_{1}... g_{n} \pd ^{-1}, g_{1}, ... , g_{n-1} ) .\]
\ \ \ The following is a particular case of \ref{lreduction} in the next sub-section:

\begin{num}
{\bf Lemma}: If $\A_{\dt}$ is a cyclic \G-sheaf such that any $\A_{n}$ is c-soft, then $HH_{*}(\G;\A_{\dt}),  HC_{*}(\G;\A_{\dt}),$ $HP_{*}(\G;\A_{\dt}) $ are computed by the cyclic vector space $\C (\G;\A_{\dt})$.
\end{num}


\subsection{Cyclic Groupoids ; Gysin sequences}

\begin{num}
\label{cycat}
{\bf Definition}: We call cyclic category an \'etale category \G\ endowed with an action of the cyclic group $\mathbb{Z}$; by this we mean there is given a continuous map $\thp:\nG{0} \rmap \nG{1} , \ c \mapsto \theta_{c}$ such that:
\begin{enumerate}
\item $\thp_{c} \in Aut \ps c \pd $, for all $c \in \nG{0}$;
\item $g\thp_{c}=\thp_{d} g$, for all $g:c \rightarrow d$ in \G.
\end{enumerate}
A morphism between two cyclic categories $(\G, \thp), (\Ha, \tau)$ is a continuous functor $f:\G \rmap \Ha$ such that $f(\thp_{c})=\tau_{f\ps c \pd} \ \ \ \forall \ \ c \in \nG{0}$.
\end{num}
\ \ \ For discrete groupoids this agrees with the old definition given by Burghelea (\cite{Bu}, page 358). There is an action of $\mathbb{Z}$ on the space of arrows: the generator acts as $\nG{1} \rmap \nG{1} \ , \ g \mapsto \thp_{t \ps g \pd}g$. The localization of $(\G,\thp)$, denoted $\G_{\ps \thp \pd}$, is obtained from \G\ by imposing the relations $\thp_{c}=id_{c}  \ \ \ ,\forall \ c \in \nG{0}$; to be more precise about the topology, put $\nG{0}_{\ps \thp \pd}:=\nG{0} \ , \ \nG{1}_{\ps \thp \pd}:=\nG{1} / \mathbb{Z}$ with the obvious structure maps. It is not difficult to see that \gth\ is still an \'etale category.
\par We call $(\G,\thp)$ elliptic if $ord(\thp_{c}) < \infty\ \ \ ,\forall \ c \in \nG{0}$. We call it hyperbolic if for any $g: c \rmap d$ in $\G$, the equality $g\thp_{c}^{n}= g$ holds just for $n= 0$ (in particular $ord(\thp_{c}) < \infty\ \ \ ,\forall \ c \in \nG{0}$; if $\G$ is a groupoid, this is the only condition). 
\par If $(\G,\thp) , \ (\G^{'},\thp^{'})$ are cyclic categories, so is $(\G \times \G^{'} ,  \thp \times \thp^{'})$. If $(\G,\thp)$ is a cyclic category   then so is $(\G, \thp^{-1})$ , provided that the map $\nG{0} \rmap \nG{1}, \ c \mapsto \thp_{c}^{-1}$ is continuous; in this case, for any other cyclic category $(\G^{'},\thp^{'})$, the localization of  $(\G,\thp^{-1}) \times (\G^{'},\thp^{'})$ is denoted by $\G \wedge \G^{'}$.

\begin{num}
\label{excyclic}
\emph{
{\bf Examples}: 
\begin{enumerate}
\item if \G\ is an \'etale category, then $(\G, id)$ is a cyclic category and $\G_{\ps id \pd}=\G$.
\item $(\Lambda_{\infty},T)$, where $T([n]):=t_{n}^{n+1}$, is a hyperbolic cyclic category with $(\Lambda_{\infty})_{\ps T \pd}=\Lambda$.
\item If $G$ is a group, $g \in Center(G)$ then $(G,g)$ is a cyclic category with $G_{\ps g \pd}=G/<g>$.
\item If $C_{r}$ is a cyclic group generated by $\ga$, $ord(\ga)=r+1$ then $(C_{r})_{\ps \ga \pd}=$trivial and $\Lambda_{\infty} \wedge C_{r}=\Lambda_{r}  \ \ \ ,\forall\ \ 0 \leq r \leq \infty$.
\item For any cyclic category $(\G,\thp)$: $\G \wedge *=\gth \ , \ \G \wedge \mathbb{Z}=\G$.
\end{enumerate}
}
\end{num}

\begin{num}
\label{le}
{\bf Lemma}: Let $(\G,\thp)$ be a hyperbolic cyclic category, $\f:\G \rmap \gth$ the projection functor. Then for any $c \in \nG{0}$, the functor:
\[ \mathbb{Z}=<\thp_{c}> \ \tilde{\rmap} c/\f  ,\]
which sends the single object $*$ of $\mathbb{Z}$ to $\ps id_{c}, c\pd$ and $n \in \mathbb{Z}$ to $\thp^{n}:\ps id_{c}, c\pd \rmap \ps id_{c}, c\pd$ is a strong deformation retract of $d/\f$. Moreover, it is a Morita equivalence if $\G$ is a groupoid.
\end{num}

\emph{proof}: Choose a set-theoretic map $\sigma : \nG{1}_{\ps \thp \pd} \rmap \nG{1}$ which is a retract of $\f$ and:
\[ s(\sigma\ps h \pd)=s\ps h\pd, \ \  t(\sigma\ps h \pd)=t\ps h\pd, \ \ \sigma \ps c\pd= c, \ \ \ \forall \ h\in \nG{1}_{\ps \thp \pd}, c \in \nG{0} .\]
We get a map $n: \G \rmap \mathbb{C}$ uniquely determined by the equality:
\[   g=\thp^{n\ps g\pd} \sigma(\f \ps g\pd), \ \ \ \forall\ g \in \nG{1} .\]
Recall (\ref{comma}) that the discrete category $c/\f$ has as objects pairs $\ps h, d\pd$ with $d \in \nG{0}, h: c \rmap d$ a morphism in $\gth$ and as morphisms from  $\ps h, d\pd$ to $\ps h\ ', d\ '\pd$ those morphisms $g: d \rmap d\ '$ in $\G$ with the property $\f \ps g \pd h= h\ '$. Define the retraction $\psi: d/\f \rmap \mathbb{Z}$ by sending a morphism $g: \ps h, d\pd \rmap \ps h\ ', d\ '\pd$ (in $c/\f$) to $n(g\sigma \ps h\pd)$. The deformation retraction $F: \psi\ \compose \f \rmap Id$ is defined as follows: to an object $\ps h, d\pd$ in $c/\f$ it associates the morphism $F\ps h, d\pd = \sigma\ps h \pd : \psi(\f\ps h, d\pd)= \ps id_{c}, c\pd \rmap \ps h, d\pd$ in $c/\f$. That $\f$ is an essential equivalence when $\G$ is a groupoid is obvious. 

\begin{prop}
\label{gysin}
: If $(\G,\thp)$ is a hyperbolic cyclic category, then for any $\A \in Sh(\gth)$ there is a long exact sequence:
\[ . \ . \ . \ \rmap H_{n}(\G;\A) \rmap H_{n}(\gth;\A) \stackrel{d}{\rmap} H_{n-2}(\gth;\A) \rmap H_{n-1}(\G;\A) \rmap .\ .\ .\ \ \ .\]
\ \ \ Here the boundary is of type $d=- \cap e(\G,\thp)$, the cap product by some cohomology class $e(\G,\thp) \in H^{2}(\gth;\mathbb{C})$ which does not depend on $\A$ and is called the Euler class of $(\G, \thp)$. Moreover, $e(\G, \thp)$ has the naturality property in the following sense: for any morphism $f: (\G, \thp) \rmap (\Ha, \tau)$ of hyperbolic cyclic categories, $f^{*}e(\Ha, \tau)= e(\G, \thp)$. 
\end{prop}

\emph{proof}: The spectral sequence for \f\ (\ref{HS}, \ref{Etcat}), $E_{p,q}^{2}=H_{p}(\gth;L_{q}\f _{\, !}\A) \Longrightarrow H_{p+q}(\G;\A)$, has:
\[ (L_{q}\f _{\, !}\A)_{c}=H_{q}(c/\f ;\A)=H_{q}(\mathbb{Z};\A_{c})\ .\]
The last equality follows from the previous lemma and \ref{inv} (or Morita invariance \ref{mor} if \G\ is a groupoid) and is continuous with respect to $c \in \nG{0}$. So we get $L_{q}\f _{\, !}\A=\A$ if $q \in \{ 0, 1\}$, and $0$ otherwise, and this implies the long exact sequence. 
\par From \ref{spseq}, \ref{HS} we know its boundary is of type $d= - \cap e(\A)$ for some $e(\A) \in Ext_{\G_{\ps \thp \pd}}(\A, \A)$. Recall (see \ref{coh}, \ref{leseq}) that the action of $H^{2}(\G_{\ps \thp \pd}; \mathbb{C})= Ext_{\G_{\ps \thp \pd}}(\mathbb{C}, \mathbb{C})$ on $H_{*}(\G_{\ps \thp \pd}; \A)$ is defined  using the morphism: $Ext_{\G_{\ps \thp \pd}}(\mathbb{C}, \mathbb{C}) \rmap Ext_{\G_{\ps \thp \pd}}(\A, \A), u \mapsto u \otimes \A\, $. So it is enough to prove that $e(\A)= e(\mathbb{C}) \otimes \A\, $.
Recall also that the spectral sequence we used is obtained from the equality: 
$ H_{*}(\G; \A)=\mathbb{H}_{*}(\gth; \el_{\dt}(\A))$, where $\el_{\dt}(\A))=\el \f_{\, !}(\A)$ is a chain complex in $Sh(\gth)$ which can be described as follows (we spell out the general definition of $\el \f_{\, !}$, \cite{CrMo}). Define:
\[X_{n}=\{ (g_{1}, . . . , g_{n+1}) : (g_{1}, . . . , g_{n}) \in \nG{n}, g_{n+1} \in \nG{1}_{\ps \thp \pd}, s(g_{n})= t(g_{n+1}) \} ,\]
and the maps $\alpha, \beta: X_{n} \rmap \nG{0} , \alpha (g_{1}, . . . , g_{n+1})=t(g_{1}), \beta(g_{1}, . . . , g_{n+1})=s(g_{n+1})$. Then $\el_{n}(\A)=\beta_{ !}\alpha^{*} \A$ has the stalk at $c \in \nG{0}$:
\[ (\el_{n}(\A))_{c}= \bigoplus_{d \in \nG{0}} \A_{d} \otimes \mathbb{C}[ \alpha^{-1}(d) \cap \beta^{-1}(c)] ,\]
and the action of $\gth$ and the boundaries are:
\[ (a, g_{1}, . . . , g_{n+1}) g = (a, g_{1}, . . . , g_{n+1}g) ,\]
\[ d_{n}(a, g_{1}, . . . , g_{n+1})= (ag_{1}, . . . , g_{n+1}) + \sum_{i=1}^{n} (-1)^{i} (a, g_{1}, . . ., g_{i}g_{i+1}, . . . , g_{n+1}) .\]
We get the following representative for $e(\A)$:
\[ 0 \rmap \A \stackrel{j}{\rmap} Coker(d_{1}) \stackrel{d_{0}}{\rmap} \el_{0}(\A) \stackrel{\epsilon}{\rmap} \A \rmap 0 ,\]
where $\epsilon \ps a, g_{1} \pd = ag_{1}, j \ps a \pd = \widehat{\ps a, \thp, 1 \pd}$. Denoting this extension by $u(\A)$, remark there is a map:
\[ \el_{\dt}(\A) \rmap \el_{\dt}(\mathbb{C}) \otimes \A, \ \ (a, g_{1}, . . . , g_{n+1}) \mapsto (ag_{1} . . . g_{n+1}, g_{1}, . . . , g_{n+1}) \]
which induces a map $u(\A) \rmap u(\mathbb{C}) \otimes \A$. This proves $e(\A)=e(\mathbb{C}) \otimes \A$.
\par Denote by $\tilde{\el}_{\dt}(\mathbb{C}), \tilde{u}(\mathbb{C})$ the analogous constructions for $(\Ha, \tau)$ and define a chain-map:
\[ \el_{\dt}(\mathbb{C}) \rmap f^{*}\tilde{\el}_{\dt}(\mathbb{C}), \ \ (g_{1}, . . . , g_{n+1}) \mapsto (f\ps g_{1} \pd, . . . , f\ps g_{n+1} \pd) .\]
It induces a map of extensions $u(\mathbb{C}) \rmap f^{*}\tilde{u}(\mathbb{C})$ and this proves the naturality.

\begin{num}
\emph{
{\bf Remark}: Applying \ref{gysin} to $(\Lambda_{\infty},T)$ \ (see \ref{excyclic}) we get the usual $SBI$-sequence for cyclic vector-spaces as the Gysin-sequence for the projection $\Lambda_{\infty} \rmap \Lambda$ and the $S$-boundary as an Euler class $e(\Lambda_{\infty},T) \in H^{2}(\Lambda; \mathbb{C})$. Our description of $S$ in terms of extensions (as described in the previous proof) is very close to the one given in \cite{Ni}, pp. 565.
}
\end{num}

\begin{num}
\label{twist}
{\bf Definition}: Let $(\G,\thp)$ be a cyclic category. A $\thp$-cyclic \G-sheaf is a $\infty$-cyclic object $\A_{\dt}$ in $Sh(\G)$ (i.e. a contravariant functor $\Lambda_{\infty} \rmap Sh(\G)$ \ cf \ref{ciclic}) such that, for any $c \in \nG{0}$, the morphism $(t_{c, n})^{n+1}:(\A_{n})_{c} \rmap (\A_{n})_{c}$ coincides with the action of $\thp_{c}$. In other words, a \thp-cyclic object is a $\Lambda_{\infty} \wedge \G$-sheaf. Define:
\[ HH_{*}(\G,\thp;\A_{\dt}):=H_{*}(\Lambda_{\infty} \times \G;\A_{\dt}) \ , \ HC_{*}(\G,\thp;\A_{\dt}):=H_{*}(\Lambda_{\infty} \wedge \G;\A_{\dt}) .\]
\end{num}
\ \ \ If \G\ is an \'etale groupoid, $\thp=id$, this agrees with the earlier definition of cyclic \G-sheaves (\ref{ciclicsheaves}) and their homologies (\ref{nustiu}).

\begin{num}
\emph{
{\bf $SBI$-sequence}: Since $\Lambda_{\infty} \times \G$ is always hyperbolic, \ref{gysin} applies, so for any \thp-cyclic \G-sheaf $\A_{\dt}$ there is a long exact sequence: 
\[ .\ .\ .\  \rmap HH_{n}(\G ,\thp ;\A_{\dt}) \rmap HC_{n}(\G ,\thp ;\A_{\dt}) \rmap HC_{n-2}(\G ,\thp ;\A_{\dt}) \rmap HH_{n-1}(\G ,\thp ;\A_{\dt})  \rmap .\ .\ .\ \]
}
\end{num}

\begin{num}
\label{important}
\emph{
{\bf Remark}: If $(\G,\thp)$ is a cyclic category and  $\A_{\dt}$ is a $\thp$-cyclic \G-sheaf, then $HH_{*}(\G,\thp;\A_{\dt})=\mathbb{H}_{*}(\G;(\A_{\dt},b))$ (one way to see this is by applying the spectral sequence for the projection $\Lambda_{\infty} \times \G \rmap \G$). This implies that the analogue of \ref{SBIARG} holds. 
}\
\end{num}

\begin{num}
\label{draci}
{\bf Lemma}: If $(\G,\thp)$ is a cyclic category and  $\A_{\dt}$ is a $\thp$-cyclic \G-sheaf, then  the $\infty$-cyclic vector space $n \mapsto H_{q}(\G;\A_{n})$ is cyclic for every $q$, and  there are spectral sequences:
\[ E_{p,q}^{2}=HH_{p}(H_{q}(\G;\A_{\dt})) \Longrightarrow HH_{p+q}(\G,\thp;\A_{\dt}) \ , \ E_{p,q}^{2}=HC_{p}(H_{q}(\G;\A_{\dt})) \Longrightarrow HC_{p+q}(\G,\thp;\A_{\dt}) .\]
\end{num}

\emph{proof}: Consider the diagram with columns coming from cyclic categories:
\[ \xymatrix{
\Lambda_{\infty} \times *  = & \Lambda_{\infty} \ar[d]_{\pi} & \Lambda_{\infty} \times \G \ar[l]_{\pi_{1}} \ar[r]^{\pi_{2}} \ar[d]_{\psi} & \G \ar[d]_{\f} &  =  * \times \G \\
\Lambda_{\infty} \wedge *  = & \Lambda & \Lambda_{\infty} \wedge \G \ar[l]_{\pi_{1}^{\, '}} \ar[r]^{\pi_{2}^{\, '}} & \gth & =  * \wedge \G
}
\]
Here $\pi_{1}, \pi_{2}, \pi_{1}^{\, '}, \pi_{2}^{\, '}$ are the projections. For any integer $n \geq 0$ we get from \ref{le} strong deformation retracts: 
\[ [n]/\pi \tilde{\leftarrow} \mathbb{Z} , \  \ [n]/\pi_{1} =([n]/\Lambda_{\infty}) \times \G \tilde{\leftarrow} \G , \ \ [n]/\pi_{1}^{\,'} = ([n]/\pi) \wedge \G \tilde{\leftarrow} \mathbb{Z} \wedge \G=\G .\]
 The spectral sequences induced by $\pi_{1}$ and $\pi_{1}^{\, '}$ together with \ref{inv}, give the desired spectral sequences.

\begin{num}
\label{reduction}
\emph{ Assume that $(\G,\thp)$ is a cyclic groupoid and $\A_{\dt}$ is a $\thp$-cyclic \G-sheaf such that each $\A_{n}$ is c-soft. From \ref{important} we see that $HH_{*}(\G,\thp;\A_{\dt})$ is computed by the bi-simplicial vector space $B_{\dt}(\G;\A_{\dt})$; so, as in \ref{atunits}, it is computed by the simplicial vector space $\C (\G ; \A_{\dt})$. Define $\C (\G, \thp ; \A_{\dt})$ as the cyclic vector space having $\C (\G ; \A_{\dt})$ as underlying simplicial vector space and the following cyclic structure:
\[ t ( a \ | \ g_{1},. . ., g_{n} ) = ( t\ps a \pd \thp^{-1}g_{1} ... g_{n} | \ps g_{1}... g_{n} \pd ^{-1}\thp , g_{1}, ... , g_{n-1} ) .\]
}
\end{num}

\begin{num}
\label{lreduction}
{\bf Lemma and definition}: If $(\G,\thp)$ is a cyclic groupoid and $\A_{\dt}$ is a \thp-cyclic \G-sheaf such that each $\A_{n}$ is c-soft, then $HH_{*}(\G,\thp ;\A_{\dt}),  HC_{*}(\G,\thp ;\A_{\dt})$ are computed by the cyclic vector space $\C (\G, \thp ;\A_{\dt})$. In this case, define $HP_{*}(\G,\thp ;\A_{\dt})=HP_{*}(\C (\G, \thp ;\A_{\dt}))$.  
\end{num}

\emph{proof}: On the standard resolution of $\mathbb{C}$ in $Sh(\G)$ (see \ref{standard}), considered as a simplicial \G-sheaf, we define a structure of $\thp$-cyclic \G-sheaf by:
\[ t ( g_{0},. . ., g_{n} ) =(\thp^{-1}g_{0}g_{1}... g_{n}, \ps g_{1}... g_{n} \pd^{-1}\thp , g_{1}, . . . , g_{n-1}) .\] 
Denote it by $\B_{\dt} \in Sh(\Lambda_{\infty} \wedge \G)$, and define $\el_{\dt} \in Sh(\Lambda_{\infty} \wedge \G)$ by $\el_{n}=\A_{n} \otimes \B_{n}$ with the structure maps given by the tensor-product of the structure maps of $\B_{\dt}$ and $\A_{\dt}$. It is a standard fact that the projection $\el_{\dt} \rmap \A_{\dt}$ induces a quasi-isomorphism $(\el_{\dt}, b) \simeq (\A_{\dt}, b)$ (it follows, for instance, from the fact that $\B_{\dt}$ is a resolution of $\mathbb{C}$ and from the Eilenberg-Zilber theorem); the "$SBI$-trick" (see \ref{important}) implies that $HC_{*}(\G,\thp;\el_{\dt}) \stackrel{\sim}{\rmap} HC_{*}(\G,\thp;\A_{\dt})$. From this and \ref{draci} applied to the $\thp$-cyclic $\G$-sheaf $\el_{\dt}$, we get a spectral sequence with: $E_{p,q}^{2}=HC_{p}(H_{q}(\G;\el_{\dt})) \Longrightarrow HC_{p+q}(\G,\thp;\A_{\dt}).$
It is enough to make a straightforward remark: $H_{q}(\G;\el_{n})=0,$ for all $q \geq 1$ (all $\el_{n}$'s are free, \cite{CrMo}) and for $q=0$, the cyclic vector space $n \mapsto H_{0}(\G;\el_{n})$ is in fact $\C (\G, \thp ;\A_{\dt})$.

\begin{prop}
\label{celliptic}
(elliptic case): Let $(\G,\thp)$ be an elliptic cyclic groupoid, $\A_{\dt}$ a \thp-cyclic \G-sheaf. For any integer $n$ put $\A_{\ps \thp \pd, n}:=\f _{\, !}(\A_{n}) \in Sh(\gth)$, where $\f:\G \rmap \gth$ is the projection. Then $\A_{\ps \thp \pd, \dt}$ is a cyclic $\gth$-sheaf and:
\[ HH_{*}(\G,\thp ;\A_{\dt})=HH_{*}(\gth ;\A_{\ps \thp \pd, \dt}) \ \ , \ \ HC_{*}(\G,\thp ;\A_{\dt})=HC_{*}(\gth ;\A_{\ps \thp \pd, \dt }) ,\] and the analogue for $HP_{*}$.
\end{prop}

\emph{proof}: Assume for simplicity that each $\A_{n}$ is c-soft (in general we work with c-soft resolutions). From \ref{le} we have the Morita equivalences $<\thp_{c}> \ \stackrel{\sim}{\rmap} c/\f \ \ ,\ \forall c \in \nG{0}$; \ref{HS} gives $(L_{q}\f_{\,!}\A_{n})_{c}=$ \\ $H_{q}(<\thp_{c}>;(\A_{n})_{c})=0 \ ,\ \forall q \neq 0$ and $(\A_{\ps \thp \pd, n})_{c}=Coinv_{\thp_{c}}((\A_{n})_{c})$. In particular, the $\infty$-cyclic object $\A_{\ps \thp \pd, \dt}$ is cyclic. Also the spectral sequence of \f\ degenerates; this ensures that the obvious projection of bi-simplicial vector spaces $\{B_{p}(\G;\A_{q}) : p, q \geq 0 \} \rmap \{ B_{p}(\gth ; \A_{\ps \thp \pd , q}) : p, q \geq 0 \}$  is a quasi-isomorphism on the $q=constant$ columns. By the Eilenberg-Zilber theorem, it is a quasi-isomorphism between their diagonals, i.e. the projection $\C (\G, \thp ; \A_{\dt}) \rmap \C (\gth, id;\A_{\ps \thp \pd, \dt})$ of cyclic vector spaces induces isomorphism in their Hochschild homologies; now \ref{lreduction} ends the proof.  
\begin{prop}
\label{chyperbolic}
(hyperbolic case): Let $(\G,\thp)$ be a hyperbolic cyclic groupoid, $\A_{\dt}$ a \thp-cyclic \G-sheaf.Then $\widetilde{HH}_{*} (\A_{\dt})\in Sh(\G)$ are in fact \gth-sheaves, and there are spectral sequences:
\[  E_{p,q}^{2}=H_{p}(\G;\widetilde{HH}_{*}(\A_{\dt}))   \Longrightarrow HH_{*}(\G,\thp;\A_{\dt}) \ \ , \ \
    E_{p,q}^{2}=H_{p}(\gth;\widetilde{HH}_{*}(\A_{\dt}))   \Longrightarrow HC_{*}(\G,\thp;\A_{\dt}) .\]
\ \ \ Moreover, $HC_{*}(\G,\thp;\A_{\dt})$ are modules over the ring $H^{*}(\gth;\mathbb{C})$ and $S$ in the $SBI$ sequence is the (cap-) product by the Euler class $e(\G, \thp) \in H^{2}(\gth;\mathbb{C})$.
\end{prop}

\emph{proof}: That \thp\ acts trivially on $\widetilde{HH}_{*} (\A_{\dt})$ follows by using the homotopy between $id$ and \thp\ , $h:\A_{\ \dt} \rmap \A_{\ \dt+1}$ given by:
\[ h_{n}=s_{-1} \ps 1+\tau + . . . + \tau^{n} \pd : \A_{n} \rmap \A_{n+1}\ .\]  
From the construction in the proof of \ref{le}, we get strong deformation retracts:
\[ c/\pi_{2}=\Lambda_{\infty} \times (c/\G) \tilde{\leftarrow} \Lambda_{\infty} , \ \ \ c/\pi_{2}^{\, '}=\Lambda_{\infty} \wedge (c/\f) \tilde{\leftarrow} \Lambda_{\infty} \wedge <\thp_{c}>=\Lambda_{\infty} .\]
Now the spectral sequences follow from the spectral sequences \ref{HS} induced by $\pi_{2}, \pi_{2}^{\, '}$ (see the diagram in the proof of \ref{draci}), and \ref{inv}. The last part follows from \ref{gysin} applied to $\Lambda_{\infty} \times \G$ and the isomorphism $H^{*}(\Lambda_{\infty} \wedge \G;\mathbb{C}) \simeq H^{*}(\gth ;\mathbb{C})$. The last isomorphism follows from the spectral sequence induced by $\pi_{2}^{\, '}$ in cohomology (cf. \cite{Ha3,Mo4,SGA}).  

\begin{num}
\emph{
{\bf Eilenberg-Zilber Theorem in the cyclic case}: As an application we give a new proof of the spectral sequence for cylindrical objects which is one of the main results in \cite{GeJo}. Recall (\cite{GeJo}, pp. 164) that a cylindrical vector space is a functor $C_{\dt \dt}:\Lambda_{\infty} \wedge  \Lambda_{\infty} \rmap \underline{Vs}$. Its diagonal is naturally a cyclic vector space. }
\end{num}
\par
\begin{num}
{\bf Corollary}: If $C_{\dt \dt}$ is a bi-cyclic vector space which is cylindrical, then the $\infty$-cyclic vector space $n \mapsto H_{q}(C_{n, \dt})$ is cyclic for any $q \geq 0$,  and there are spectral sequences:
\[ E_{p,q}^{2}=HH_{p}(H_{q}(C_{\dt \dt})) \Longrightarrow HH_{p+q}(diag(C_{\dt \dt })) ,\]
\[ E_{p,q}^{2}=HC_{p}(H_{q}(C_{\dt \dt})) \Longrightarrow HC_{p+q}(diag(C_{\dt \dt}))  .\]
\end{num}

\emph{proof}: With the same method as in the proof of \ref{lreduction} we see that the homologies $HH_{*}(\Lambda_{\infty},T;C_{\dt \dt})$,  $HC_{*}(\Lambda_{\infty},T;C_{\dt \dt})$ are computed by the cyclic module $diag(C_{\dt \dt})$. Then apply \ref{draci} to $(\Lambda_{\infty},T)$.


\subsection{Cyclic Homology of crossed products by \'etale groupoids}

\begin{num}
\label{crpr}
\emph{ Recall that if $G$ is a discrete group acting on an unital algebra $A$, the cross-product algebra $A \cross G$ is $A \otimes \mathbb{C}[G]$ with the convolution product $\ps a, g \pd \ps b, h \pd = (\ps ah \pd b, gh)$. We see that the induced cyclic vector space $(A \cross G)^{\natural}$  (see \ref{alfaex}) has:
}
\end{num}
\[  (A \cross G)^{\natural}\ps n \pd = A^{\otimes \ps n+1 \pd } \otimes \mathbb{C}[G^{n+1}] ,\]
while the cyclic structure is given by the formulas (see also \cite{Ni1}):
\[ d_{i}( a_{0}, . . . , a_{n} | g_{0} , . . . , g_{n} ) = \left\{ \begin{array}{ll}
           ( a_{0}, . . . , \ps a_{i}g_{i+1} \pd a_{i+1} , . . . , a_{n} | g_{0} , . . . , g_{i}g_{i+1} , . . . , g_{n} ) & \mbox{if $0 \leq i \leq n-1$} \\
           ( \ps a_{n}g_{0} \pd a_{0} , a_{1} , . . . , a_{n-1} | g_{n}g_{0} , g_{1} , . . . , g_{n-1} ) & \mbox{if $i=n$}
                                                                       \end{array}
           \right. ,\]
\[ s_{i}( a_{0}, . . . , a_{n} | g_{0} , . . . , g_{n} ) =( . . . , a_{i} , 1 , a_{i+1} , . . . | . . . , g_{i} , 1 , g_{i+1} , . . . ) ,\]
\[ t( a_{0}, . . . , a_{n} | g_{0} , . . . , g_{n} ) = ( a_{n} , a_{0}, . . . , a_{n-1} | g_{n} , g_{0} , . . . , g_{n-1} ) .\]

\begin{num}
{\bf Lemma and definition}: Let $\G$ be an \'etale groupoid such that $\nG{0}$ is Hausdorff, $\A$ a c-soft $\G$-sheaf of complex algebras. If $u, v \in \Gamma_{c}(\G; s^{*}\A)$, then the following formula:
\[ (u * v)\ps g \pd= \sum_{g_{1}g_{2}=g} (u \ps g_{1}\pd g_{2})v\ps g_{2}\pd  \ \ \forall g \in \G ,\]
gives a well defined element $u * v \in \Gamma_{c}(\G; s^{*}\A)$. The resulting algebra, $(\Gamma_{c}(\G; s^{*}\A), *)$ is called the crossed product of $\A$ and $\G$ and is denoted $\A\cross_{alg}\G$.
\end{num}

\emph{proof}: The element $u * v \in \Gamma_{c}(\G; s^{*}\A)$ is obtained via the composition of:
\[ \Gamma_{c}(\G; s^{*}\A) \times \Gamma_{c}(\G; s^{*}\A) \rmap \Gamma_{c}(\G \times \G; s^{*}\A \boxtimes s^{*}\A) \rmap \Gamma_{c}(\G \times_{\nG{0}} \G; s^{*}\A \boxtimes s^{*}\A\, |_{\G \times_{\nG{0}} \G}) \rmap \Gamma_{c}(\G; s^{*}\A) \]
where the first map is the obvious one, the second is the restriction to $\G \times_{\nG{0}} \G$ (which is closed in $\G \times \G$ since $\nG{0}$ is Hausdorff), and the third is:
\[ (a,b \ | \ g_{1}, g_{2}) \mapsto (\ps ag_{2}\pd b, g_{1}g_{2}). \]

\begin{num}
\label{coco}
\emph{
{\bf Examples}:
\begin{enumerate}
\item If $G$ is a discrete group acting on an algebra $A$, then $A\cross_{alg}G$ is the usual crossed product (described in \ref{crpr}).
\item If $\G$ is a smooth \'etale groupoid, $\A=\C_{\nG{0}}^{\infty}$, then $\A\cross_{alg}\G=\conv$ is the convolution algebra of $\G$ (we use that $\conv=\Gamma_{c}(\G; \C_{\G}^{\infty})$ and $ \C_{\G}^{\infty}= s^{*}\C_{\nG{0}}^{\infty}$ since $s$ is \'etale). We will come back to this example in the next subsection.
\end{enumerate}
}
\end{num}

\begin{num}
\label{cycrpr}
{\bf Definition}: Let \G\ be an \'etale groupoid, \A\ a c-soft \G-sheaf of complex unital algebras.Consider Burghelea's space (\cite{BrNi}):
\[ B^{\ps n \pd }=\{ \ps \ga_{0} , \ga_{1} , . . . , \ga_{n} \pd \in \nG{n+1} : t\ps \ga_{0} \pd = s\ps \ga_{n} \pd \} ,\]
and $\sigma_{n}:B^{\ps n \pd } \rmap (\nG{0})^{n+1} , \sigma_{n} ( \ga_{0} , . . . , \ga_{n} ) = ( s \ps \ga_{0} \pd , . . . , s \ps \ga_{n} \pd) $. Define the cyclic vector space $(\A \cross \G)^{\natural}$  by:
\[ (\A \cross \G)^{\natural}(n) = \gc(B^{\ps n \pd } ; \sigma_{n}^{*}\A^{\boxtimes \ps n+1 \pd }) .\]
\ \ \ For the structure maps, we keep the same formulas as in  \ref{crpr}. Its homologies are denoted $HH_{*}(\A \cross \G) , HC_{*}(\A \cross \G) , HP_{*}(\A \cross \G).$ 
\end{num}

\ \ \ Maybe this definition  is not very transparent. The main examples are the discrete case \ref{crpr} and the smooth case (see Proposition \ref{redloops}). A way to think of $(\A \cross \G)^{\natural}$ is as follows (and this is in fact the main motivation). Usually our objects are endowed with topologies; in particular, the relevant cyclic homology of $\A\cross_{alg}\G$ uses the topology (i.e. the projective tensor product \cite{Gr} is used for defining $(\A\cross_{alg}\G)^{\natural}$). We have inclusions:
\[ (\A \cross_{alg} \G)^{\otimes \ps n+1\pd } \hookrightarrow \gc (\G^{n+1}; s_{n+1}^{*}(\A^{\boxtimes \ps n+1\pd}))  \ \ , \ (s_{n+1}=s \times . . . \times s \ ; \ps n+1 \pd \ times) \] 
and the vector spaces  $\gc (\G^{n+1}; s^{*}_{n+1}(\A^{\boxtimes \ps n+1\pd}))$ are the best candidates for defining $(\A \cross_{alg} \G)^{\natural}$ by taking into account the topology. From the first half of \ref{basicsgc}.2 we get (compare to Proposition 4.1 in \cite{BrNi}):

\begin{num}
\label{lele}
{\bf Lemma}: $(\A \cross \G)^{\natural}(n)=\gc (\G^{n+1}; s^{*}_{n+1}(\A^{\boxtimes \ps n+1\pd} ))/ \{ u : u\,|_{B^{\ps n \pd}} = 0 \}.$
\end{num}
\par So our definition can be viewed as a ``topological normalization'' of the cyclic vector space associated to $\A \cross _{alg} \G $. Of course, in practice one more step is needed: to prove that passing to this ``topological normalization'' does not change the cyclic homology (see for instance \ref{redloops}).
\par What is also very important in our choice is that $(\A\cross\G)^{\natural}$ is the maximal object with the property that the formulas in \ref{crpr} make sense; this allows us to use our methods in this general setting. And as a final motivation for our choice and notation, let us just remark that the usual formulas for the Chern-character (see e.g. \cite{Co3,Lo}) define a Chern-Connes character: $Ch: K_{*}^{alg}(\A \cross_{alg} \G) \rmap HP_{*}(\A \cross \G)$. \par Because of these we do believe that the results (or better: the methods) we describe here may be useful to a larger extent.

\begin{num}
\label{NNN}
\emph{
{\bf The cyclic groupoid \Z; the sheaves $\A^{\natural}_{\nB{0}}$}: An important role in understanding $(\A \cross \G)^{\natural}$ (even in the discrete case \ref{crpr}) is played by the \G-space of loops $\nB{0}=\{ \ga \in \nG{1} : s \ps \ga \pd = t \ps \ga \pd \}$; the action of \G\ is given by:
\[ \nB{0} \times_{\nG{0}} \G \rmap \nB{0} \ , \ \ps \ga , g \pd \mapsto g^{-1} \ga g .\]
Let $\Z:=\nB{0} \cross \G$ be the associated groupoid (see \ref{actions}); it is a cyclic groupoid with the cyclic structure $\thp$ defined by $\thp (\ga ):=(\ga , \ga)$. Define a $\thp$-cyclic \Z-sheaf $\A^{\natural}_{\nB{0}}$ by $\A^{\natural}_{\nB{0}}:=s^{*}(\A^{\natural})$ where $s:\nB{0} \rmap \nG{0}$ is the restriction of the source map, and $\A^{\natural}=\A^{\otimes \ps *+1 \pd}$ was defined in \ref{standard}; the formulas for the structure maps are defined in such a way that, at the stalk at $\ga \in \nB{0}$ with $s\ps \ga \pd =c$, the $\infty$-cyclic vector space $(\A^{\natural}_{\nB{0}})_{\ga}$ is $\A^{\natural}_{c, \sigma_{\ga}}$, i.e. the one associated to the action of $\ga$ on the algebra $\A_{c}$ (see \ref{alfaex}). Compare to (9), (10) in \cite{BrNi}.
}
\end{num}

\begin{prop}
\label{redcross}
: $HH_{*}(\A \cross \G)=HH_{*}(\Z , \thp ; \A^{\natural}_{\nB{0}}) \ , \ HC_{*}(\A \cross \G)=HC_{*}(\Z , \thp ; \A^{\natural}_{\nB{0}})$,\\ $HP_{*}(\A \cross \G)=HP_{*}(\Z , \thp ; \A^{\natural}_{\nB{0}})$.
\end{prop}

\emph{proof}: The $n$'th nerve of \Z\ is:
\[ \nZ{n}= \{ \ps \ga | \ga_{1} , . . . , \ga_{n} \pd : \ \ps \ga , \ga_{1} \pd \in \nZ{1} , \ps \ga_{1} , . . . , \ga_{n} \pd \in \nG{n} \} \]
(here "$ | $" is just a notation in order to separate the loops from the usual arrows). The isomorphism of vector spaces:
\[ \gc(B^{\ps n \pd } ; \sigma_{n}^{*}\A^{\boxtimes \ps n+1 \pd }) \simeq \gc ( \nZ{n} ; \A^{\ps n\pd}_{\nB{0}}) ,\]
\[ (a_{0} , . . . , a_{n} | \ga_{0} , . . . , \ga_{n} ) \mapsto ( a_{0}\ga_{1} . . . \ga_{n}\ga_{0} , a_{1}\ga_{2} . . . \ga_{n}\ga_{0} , .\ .\ .\ , a_{n}\ga_{0} | \ga_{1}. . . \ga_{n}\ga_{0} , \ga_{1} , .\ .\ .\ ,\ga_{n} ) ,\]
(compare to 4.1 in \cite{BrNi}) gives an isomorphism  of cyclic vector spaces:
\[ (\A \cross \G)^{\natural} \simeq \C (\Z , \thp ; \A^{\natural}_{\nB{0}}) ,\]
so it is enough to use \ref{lreduction}.  
\begin{num}
\label{ZNO}
\emph{
{\bf The groupoids $\Z_{\OO}, \N_{\OO}$; the sheaves $\A^{\natural}_{\OO}$}: For $\OO \subset \nB{0}$, \G-invariant we define the cyclic groupoid $\Z_{\OO}:= \OO \cross \G$ (the restriction of $\Z$ to $\OO$). Its localization $(\Z_{\OO})_{\ps \thp \pd}$ (i.e. obtained from $\Z_{\OO}$ by imposing the relations $(\ga, 1)=(\ga, \ga) \ ,\ \forall \ \ga \in \OO$; see \ref{cycat}) is denoted by $\N_{\OO}$. These two groupoids play the role of the centralizer and normalizer of $\OO$ (see also subparagraph 4.3). We define a $\thp$-cyclic $\Z_{\OO}$-sheaf $\A^{\natural}_{\OO}$ in such a way that, 
stalkwise:
\[ (\A^{\natural}_{\OO})_{\ga}= \A^{\natural}_{c, \sigma_{\ga}}, \ \ \ \forall \ \ga \in \OO \ \ \ \ (c=s\ps \ga\pd) \]
(see \ref{alfaex}, \ref{sigma}.3). In other words, $\A^{\natural}_{\OO}$ is the restriction of $\A^{\natural}_{\nB{0}}$ to $\OO$. Define also the cyclic $\N_{\OO}$ -sheaf $\A^{\natural}_{(\OO)}:=(\A^{\natural}_{\OO})_{\ps \thp \pd}$ (see \ref{celliptic}); the stalk of $\A^{\ps n\pd}_{(\OO)}$ at $\ga \in \OO$ is $Coinv_{<\ga>}( \A_{c}^{\otimes \ps n+1\pd})$.
}
\end{num}

\begin{num}
\label{localization}
\emph{
{\bf Localization}: For $\OO$ as before, we can define the cyclic vector space $(\A \cross \G)^{\natural}_{\OO}$ as in \ref{cycrpr} with the only difference we replace $\nB{n}$ by $\nB{n}_{\OO}=\{ \ps \ga_{0} , \ga_{1} , . . . , \ga_{n} \pd \in \nB{n} : \ \ga_{0} \ga_{1}  . . .  \ga_{n} \in \OO \}$. Denote its homologies by $HH_{*}(\A \cross \G)_{\OO},HC_{*}(\A \cross \G)_{\OO}, HP_{*}(\A \cross \G)_{\OO}$. As in \ref{redcross} we have: 
\[ HH_{*}(\A \cross \G)_{\OO}= HH_{*}(\Z_{\OO} , \thp ; \A^{\natural}_{\OO})\ ,\ HC_{*}(\A \cross \G)_{\OO}= HC_{*}(\Z_{\OO} , \thp ; \A^{\natural}_{\OO})\ , \] and the analogue for $HP_{*}$. If $\OO$ is open in $\nB{0}$, there are ``extension by $0$'' maps:
$HH_{*}(\A \cross \G)_{\OO} \rmap HH_{*}(\A \cross \G);$
if $\OO$ is closed in $\nB{0}$, there are ``restriction'' maps:
$HH_{*}(\A \cross \G) \rmap HH_{*}(\A \cross \G)_{\OO}.$
The same discussion applies to $HC_{*}, HP_{*}$. The following is obvious:
}
\end{num}

\begin{prop}
 If $\nB{0}=\bigcup {\OO}$ is a \G-invariant disjoint open covering of $\nB{0}$, then:
\[ HH_{*}(\A \cross \G)=\bigoplus_{\OO} HH_{*}(\A \cross \G)_{\OO} , \ HC_{*}(\A \cross \G)=\bigoplus_{\OO} HC_{*}(\A \cross \G)_{\OO}  ,\] and the analogue for $HP_{*}$. 
\end{prop}

\begin{num}
\label{1elliptic}
\emph{
{\bf Elliptic case}: We call $\OO \subset \nB{0}$ elliptic if it is \G-invariant and $ord(\ga) < \infty$, for all  $\ga \in \OO$. From \ref{redcross},\ref{celliptic}, \ref{1st}, \ref{2nd} we get (compare to A6.2 in \cite{FeTs}, section 2 in \cite{Ni1}): 
}
\end{num}

\emph{
{\bf Theorem}(elliptic case): If $\OO$ is elliptic, then $\A^{\natural}_{(\OO)}$ is a cyclic $\N_{\OO}$-sheaf and $HH_{*}(\A \cross \G)_{\OO}=HH_{*}(\N_{\OO}; \A^{\natural}_{(\OO)}) , \ HC_{*}(\A \cross \G)_{\OO}=HC_{*}(\N_{\OO}; \A^{\natural}_{(\OO)}) ,\ HP_{*}(\A \cross \G)_{\OO}=HP_{*}(\N_{\OO}; \A^{\natural}_{(\OO)})$. In particular, there are Feigin-Tsygan-Nistor type spectral sequences:
\[ E_{p,q}^{2}=H_{p}(\N_{\OO}; \widetilde{HC}_{q}(\A^{\natural}_{(\OO)})) \Longrightarrow HC_{p+q}(\A \cross \G)_{\OO} ,\]
\[ E_{p,q}^{2}=HC_{p}(H_{q}(\N_{\OO}; \A^{\natural}_{(\OO)})) \Longrightarrow HC_{p+q}(\A \cross \G)_{\OO} ,\]
and the analogue for $HH_{*}$.
}
\newline
\par The localizations at units (i.e. at $\OO=\nG{0}$) are usually denoted by the subscript $[1]$ instead of $\OO$. We get they are equal to  $HH_{*}(\G;\A^{\natural}), \ HC_{*}(\G;\A^{\natural})$, $HP_{*}(\G;\A^{\natural})$ and the corresponding spectral sequences have at the second level: 
\[ H_{p}(\G; \widetilde{HH}_{q}(\A^{\natural})) \Longrightarrow HH_{p+q}(\A \cross \G)_{\pss 1 \pdd}\ \ , \ \ HH_{p}(H_{q}(\G; \A^{\natural}) \Longrightarrow HH_{p+q}(\A \cross \G)_{\pss 1 \pdd } \ ,\]
\[ H_{p}(\G; \widetilde{HC}_{q}(\A^{\natural})) \Longrightarrow HC_{p+q}(\A \cross \G)_{\pss 1 \pdd}\ \ , \ \ HC_{p}(H_{q}(\G; \A^{\natural}) \Longrightarrow HC_{p+q}(\A \cross \G)_{\pss 1 \pdd } \ .\]
 In the  case of crossed products by groups, this is 2.6 in \cite{Ni1}. We also know (from \ref{bdr}) the form of the boundaries $d_{p, q}^{2}$. This generalizes a similar result for crossed products by groups (see Prop. 3.2 in \cite{Ba}).

\begin{num}
\label{1hyperbolic}
\emph{
{\bf Hyperbolic case}: We call $\OO \subset \nB{0}$ hyperbolic if it is \G-invariant and  
 $ord(\ga) = \infty $, for all  $\ga \in \OO$. Denote by $e_{\OO} \in H^{2}(\N_{\OO}; \mathbb{C})$ the Euler class of the (hyperbolic) cyclic groupoid $\Z_{\OO}$. From \ref{redcross}, \ref{chyperbolic} we get (compare to A6.1 in \cite{FeTs}, 1.8 in \cite{Bu}, section 3 in \cite{Ni1}):
}
\end{num}

\emph{
{\bf Theorem}(hyperbolic case): If $\OO$ is hyperbolic, then the $HC_{*}(\A \cross \G)_{\OO}$ are modules over the ring $H^{*}(\N_{\OO};\mathbb{C})$ and $S$ in the $SBI$-sequence is the (cap-) product by the Euler class $e_{\OO} \in H^{2}(\N_{\OO}; \mathbb{C})$. Moreover, the  $\widetilde{HH}_{q}(\A^{\natural}_{\OO}) \in Sh(\Z_{\OO})$ are $\N_{\OO}$-sheaves and there are spectral sequences:
\[ E_{p,q}^{2}=H_{p}(\Z_{\OO}; \widetilde{HH}_{q}(\A^{\natural}_{\OO})) \Longrightarrow HH_{p+q}(\A \cross \G)_{\OO} ,\]
\[ E_{p,q}^{2}=H_{p}(\N_{\OO}; \widetilde{HH}_{q}(\A^{\natural}_{\OO})) \Longrightarrow HC_{p+q}(\A \cross \G)_{\OO} .\]
}


\section{Cyclic Homology of Smooth \'Etale Groupoids}


\subsection{Cyclic Homology of Smooth \'Etale Groupoids}

\ \ \ Let \G\ be a smooth \'etale groupoid. Recall that if \G\ is Hausdorff then its convolution algebra is defined as $\conv=\{ a:\nG{1} \rmap \mathbb{C} : a$ is compactly supported and smooth $\} $ with the convolution product: $(a b ) \ps \ga \pd = \sum_{\ga = \ga_{1} \ga_{2}} a\ps \ga_{1} \pd b \ps \ga_{2} \pd $. This is a locally convex algebra, which is non-unital in general. Its Hochschild and cyclic homologies are computed (cf. \ref{alfaex}) by using the $(b, b^{'})$-complex of $\conv^{\natural}$ and the projective tensor product. As remarked in \cite{BrNi} there are algebraic equalities:
\[ \conv^{\natural}\ps n \pd = \C_{c}^{\infty} (\G^{n+1}) ,\]
and the structure maps are: 
\[ d_{i}\ps a_{0}, . \ . \ . , a_{n}\pd=\left\{ \begin{array}{ll}
                                                       \ps a_{0}, . \ . \ .\ ,a_{i}a_{i+1}, . \ . \ .\ ,a_{n} \pd & \mbox{if $0 \leq i \leq n-1$} \\
                                                       \ps a_{n} a_{0}, a_{1}, . \ . \ . \ , a_{n-1} \pd & \mbox{if $i=n$}
                                                \end{array}
                                        \right. ,\]
\[ t \ps a_{0}, . \ . \ . \ , a_{n} \pd = \ps a_{n} , a_{0}, a_{1}, . \ . \ .\ , a_{n-1} \pd \  .\]
\par   In the general (i.e. non-Hausdorff) case, $\conv$ still makes sense as an algebra (see \ref{coco}.2). Rather than going into details concerning the topology, following \cite{BrNi}, we take the previous equalities as definition of the cyclic vector space $\conv^{\natural}$ (this is the relevant object for defining the Chern-character). It is important to emphasize here that, in contrast with \cite{BrNi} where Connes' definition of $\C_{c}^{\infty}(M)$ for a non-Hausdorff manifold $M$ is used, we work with our definition as given in \ref{our} (and this is essential). The homologies of the cyclic vector space  $\conv^{\natural}$ (see \ref{ciclic}) are denoted by $HH_{*}(\conv), HC_{*}(\conv), HP_{*}(\conv)$. \\
\par Denote $\C_{\nG{0}}^{\infty} = \A$; it is a \G-sheaf cf. \ref{standard}. In the definitions of $\A^{\otimes n} , \ \A^{\boxtimes n}$ we take into account the topology, as explained in \ref{topology}. The following is an extension of 3.2 in \cite{BrNi} to the non-Hausdorff case:

\begin{prop}
\label{redloops}
: For any smooth \'etale groupoid \G,  $HH_{*}(\conv), HC_{*}(\conv), HP_{*}(\conv)$ are isomorphic to $HH_{*}(\A \cross \G), HC_{*}(\A \cross \G), HP_{*}(\A \cross \G)$.
\end{prop}

\par  Since \G\ is \'etale, the elements $u \in \conv = \gc(\G; s^{*}\A)$ can be viewed as functions $\G\ \ni \ga \mapsto u\ps \ga \pd \in (s^{*}\A)_{\ga}=\A_{s\ps \ga \pd}$ and the convolution product becomes $(u * v)\ps \ga\pd  = \sum_{\ga = \ga_{1} \ga_{2} } (u\ps \ga_{1} \pd \ga_{2})v\ps \ga_{2} \pd$ (in other words, $\conv= \A\cross_{alg}\G$). We will simply denote $u * v$ by $uv$. In the same way, the elements: 
\[ u \in \C_{c}^{\infty}(\G )^{\natural} \ps n \pd = \C_{c}^{\infty}(\G ^{n+1})=\gc (\G ^{n+1};
s_{n+1}^{*}\A^{\boxtimes \ps n+1 \pd} ) \]
(here $s_{n}= s \times . . . \times s$ for $n$ times and we use the notations from \ref{topology}) can be viewed as functions: 
\[ \G^{n+1} \ni \ps \ga_{0} , . . . , \ga_{n} \pd \mapsto u \ps  \ga_{0} , . . . , \ga_{n} \pd \in ( \A \boxtimes . . . \boxtimes \A)_{(s\ps \ga_{0} \pd , . . . , s\ps \ga_{n} \pd)} .\]
\par \emph{This is the only way we are going to look at elements in} $\C_{c}^{\infty}(\G ^{n+1})$. With this it is straightforward to write the formulas for the $t$'s and the $d_{i}$'s in the general case.
\par
\emph{proof of \ref{redloops}}: From \ref{lele} we get a projection (restriction to the $\nB{n}$'s in fact) $\pi:\conv^{\natural} \rmap (\A \cross \G)^{\natural}$ which is compatible with the $t$'s and the $d_{i}$'s. It is enough to prove that $\pi$ is a quasi-isomorphism between both the $b$ and the $b'$ complexes. Denote $C_{\dt}=\conv^{\natural}, A_{\dt}=Ker(\pi)$; so $(\A \cross \G)^{\natural}=C_{\dt}/A_{\dt}$. We prove that  $((\A \cross \G)^{\natural}, b^{'}),$ $(C_{\dt}, b^{'})$ and $(A_{\dt}, b)$ are acyclic. For the two $b^{'}$ complexes this is standard. Indeed, for the first complex we use the degeneracy $s_{n}$ to get a contraction. For the second complex we use local units to get local contractions; more precisely, if $u \in C_{n}$ has $b^{'}u=0$, then $u=b^{'}v$ where $v=u \otimes \f \in C_{n+1}, \ \f \in \C_{c}^{\infty}(\nG{0})$ such that $\f=1$ around the compact $K=\{ s\ps \ga_{n} \pd : u \ps \ga_{0} , . . . , \ga_{n} \pd \neq 0 \}$ in $\nG{0}$.
\newline
\ \ \ We are left with $(A_{\dt}, b)$. We will construct two double complexes $\{ C_{p, q} ; p \geq -1, q \geq 0 \}, \{ A_{p, q} ; p \geq -1, q \geq 0 \}$ with $C_{-1, \dt}=C_{\dt}, A_{-1, \dt}=A_{\dt}$ and such that all the columns $A_{\dt, q}, q \geq 0$ and all the rows $A_{p, \dt}, p \geq 0$ are acyclic. Of course, it is enough to construct $A_{\dt,\dt}$ with these properties; the only role of $C_{\dt, \dt}$ is to facilitate the definition of $A_{\dt, \dt}$. For this we put for $p, q \geq 0$:
\[ C_{p, q} = \C_{c}^{\infty}( \G^{q+2} \times \G^{p}), \  A_{p, q}= \{ u \in C_{p, q} : u|_{\G^{\ps q+2 \pd} \times \G^{p}}=0 \} ,\]
and $C_{-1, \dt}=C_{\dt}, A_{-1, \dt}=A_{\dt}$. The boundaries are defined as follows:
\par  
 The boundaries for the row $p=constant$ of $C_{\dt, \dt}$ are defined by analogy with  $(C_{\dt}, b^{'}) \widehat{\otimes} C_{c}^{\infty}(\G^{p})$ i. e. :
\[ . \ . \ . \ \stackrel{d^{\, r}}{\rmap} C_{p, 1} \stackrel{d^{\, r}}{\rmap} C_{p, 0} , \ \ \  d^{\, r}=\sum (-1)^{i}d_{i}^{\, r} ,\]
\[ d_{i}^{\, r}(u_{0}, . . . , u_{q+1} ; v_{1}, . . . , v_{p})=(u_{0}, . . . , u_{i}u_{i+1}, . . . , u_{q+1} ; v_{1}, . . . , v_{p}), \ 0 \leq i \leq q .\]
We keep the same formulas to define the row $p=constant$ of $A_{\dt,\dt}$. To see that $A_{p, \dt}$ is acyclic for $p \geq 0$ remark that $(C_{p, \dt}/A_{p, \dt},    b^{'})$ and $(C_{p, \dt}, b^{'})$ are acyclic; this can be viewed in the same way as for $(C_{\dt}/A_{\dt}, b^{'}) $ and $(C_{\dt}, b^{'})$.
\par
The boundaries for the column $q=constant$ of $C_{\dt, \dt}$ are defined by analogy with the complex computing $HH_{*}(\conv, C_{c}^{\infty}(\G^{q+2}))$ together with an augmentation; more precisely, define them by:
\[ . \ . \ . \ \stackrel{d^{\, c}}{\rmap} C_{1, q} \stackrel{d^{\, c}}{\rmap} C_{0, q} \stackrel{d^{\, c}_{\ps -1 \pd}}{\rmap} C_{-1, q}=C_{q} , \ \ d^{\, c}=\sum_{i=0}^{p} (-1)^{i}d_{i}^{\, c} ,\]
\[ d_{i}^{\, c}(u_{0}, . . . , u_{q+1} ; v_{1}, . . . , v_{p})= \left\{ \begin{array}{lll}
         (u_{0}, . . . , u_{q}, u_{q+1}v_{1}, v_{2}, . . . , v_{p}) & \mbox{if $i=0$} \\
         (u_{0}, . . . , u_{q+1} ; v_{1} , . . , v_{i}v_{i+1}, . . . , v_{p}) & \mbox{if $0 \leq i \leq p-1$}  \ ,\\
         (v_{p}u_{0}, . . . , u_{q+1} ; v_{1}, . . . , v_{p-1}) & \mbox{if $i=p$}
    \end{array}
     \right. \]
\[ d_{\ps -1 \pd}^{\, c}(u_{0}, . . . , u_{q+1})=(u_{q+1}u_{0}, . . . , u_{q}) .\]
\par
The same formulas define the column $q=constant$ of $A_{\dt, \dt}$. To prove that $A_{\dt, q}$ is acyclic, assume for simplicity of notation that $q=0$. For any $\f \in C_{c}^{\infty}(\nG{0})$, define a degree 1 linear map of $C_{\dt, 0}$ by:
\[ h_{\f}:C_{p, 0} \rmap C_{p+1, 0} \ \ \ \ h_{\f}(u_{0}, u_{1} ; v_{1}, . . . , v_{p})=(u_{0}, \f ; \f u_{1} , v_{1}, . . . , v_{p}), \ \  p \geq 0 \]
\[ h_{\f}:C_{-1, 0} \rmap C_{0, 0} \ \ \ \ h_{\f}(u_{0})=(\f u_{0}, \f ) .\]
\par
For any $u \in C_{p, 0}=C_{c}^{\infty}(\G^{p+2})$ we have the naive formulas (correct in the Hausdorff case):
\[ (d^{c}h_{\f}u)(\ga_{0}, . . . , \ga_{p+1})=\f^{2}(t \ps \ga_{1} \pd )u(\ga_{0}, . . . , \ga_{p+1}) + \f(\ga_{1})\f (t \ps \ga_{2} \pd ) (d^{c}u) (\ga_{0}, \ga_{2}, . . . , \ga_{p+1}) ,\]
which can be written in general:
\[ (d^{c}h_{\f}u)(\ga_{0}, . . . , \ga_{p+1})=(1 \otimes germ_{t\ps \ga_{1} \pd}(\f^{2}) \otimes 1 \otimes . . . \otimes 1) * u(\ga_{0}, . . . , \ga_{p+1}) + \]
\[ + i^{\ps 2 \pd}_{germ_{\ga_{1}}\ps \f \pd}((1 \otimes germ_{t\ps \ga_{2} \pd}(\f) \otimes 1 \otimes . . . \otimes 1) * (d^{c}u) (\ga_{0}, \ga_{2}, . . . , \ga_{p+1}))   .\] 
\ \ \ Here $i_{v}^{\ps 2 \pd }$ denotes "inserting $v$ on the second place" and "$*$" is the stalkwise product. In general we do not have: $h_{\f}(A_{p, 0}) \subset A_{p+1, 0}$. Fix a metric $\rho$ defining the topology of $\nG{0}$. Take $u \in A_{p, 0}$ to be a cycle. From \ref{basicsgc} we can choose $\epsilon > 0$ such that:
\[ \rho(s\ps \ga_{0} \pd, t \ps \ga_{1} \pd ) < \epsilon  \Rightarrow u(\ga_{0}, \ga_{1} , . . . , \ga_{p})=0 .\]
\ \ \ For any $\f \in C_{c}^{\infty}(\nG{0})$ with $diam(supp \f ) < \epsilon/6$ we see that $h_{\f}(u) \in A_{p+1, 0}$. Now $u=d^{c}v$ is the boundary of $v=\sum h_{\f_{i}}u$ where $\f_{i} \in C_{c}^{\infty}(\nG{0})$ is a finite set of functions as before such that $\sum \f_{i}^{2}=1$ around the compact  $K=\{ t\ps \ga_{1} \pd : u \ps \ga_{0} , . . . , \ga_{n} \pd \neq 0 \} \subseteq \nG{0}$.  
\newline
\par The following four theorems are generalizations of the computations performed by Burghelea for groups (Theorem I in \cite{Bu}).
In the case of smooth Hausdorff \'etale groupoids, the computation of the elliptic components was done by Brylinski and Nistor: Theorem 5.6 in \cite{BrNi} computes these components in terms of homology of some double complexes. Our Theorem \ref{last} is an extension of that result, in a slightly more precise form. Emphasize also that our proof of Theorem \ref{last}, in contrast to the one given in \cite{BrNi}, will make use of the quasi-isomorphism ensured by \ref{stability} just for $\f= id$. 
\par From \ref{localization} we get (see also Proposition 3.3 in \cite{BrNi}):

\begin{st}
\label{ultimaa}
(localization): For $\OO \subset \nB{0}$ a $\G$-invariant subset, the localized homologies at $\OO$ are defined. We have linear maps $HH_{*}(\conv)_{\OO} \rmap HH_{*}(\conv)$ if $\OO$ is open in $\nB{0}$, and $HH_{*}(\conv) \rmap HH_{*}(\conv)_{\OO}$   if $\OO$ is closed in $\nB{0}$, and the same applies to $HC_{*}, HP_{*}$. 
\par If $\nB{0}=\bigcup {\OO}$ is a \G-invariant disjoint open covering of $\nB{0}$, then:
\[ HH_{*}(\conv)=\bigoplus_{\OO} HH_{*}(\conv)_{\OO} , \ HC_{*}(\conv)=\bigoplus_{\OO} HC_{*}(\conv)_{\OO} ,\]
and the analogue for $HP_{*}$; moreover, everything is compatible with the $SBI$-sequences.
\end{st}

\begin{st}
\label{ultima}
(localization at units): For any smooth \'etale groupoid \G:
\[ HP_{0}(\conv)_{\pss 1 \pdd}=\prod_{k-even} H_{k}(\G; \mathbb{C}), \]
\[ HP_{1}(\conv)_{\pss 1 \pdd}=\prod_{k-odd} H_{k}(\G; \mathbb{C}). \]
\end{st}

\emph{proof}: First of all, using \ref{stability} for $\f=id$ we get a quasi-isomorphism of sheaves on $\nG{0}$:
\[   ( \A^{\natural}, b )  \stackrel{\sim}{\rmap} (\Omega^{*}_{\nG{0}}, 0 ) . \]
It is well known (see \cite{Co2}) that it makes the $B$-boundary compatible with de Rham boundary. This, \ref{redloops}, \ref{1elliptic} for $\OO=\nG{0}$ and \ref{uit} give the proof.
\begin{st}
\label{last}
(elliptic case): For any smooth \'etale groupoid \G, and any elliptic $\OO \subset \nB{0}$:
\[ HP_{0}(\conv)_{\OO}=\prod_{k-even} H_{k}(\N_{\OO}; \mathbb{C})=\prod_{k-even} H_{k}(\Z_{\OO}; \mathbb{C}), \]
\[ HP_{1}(\conv)_{\OO}=\prod_{k-odd} H_{k}(\N_{\OO}; \mathbb{C})=\prod_{k-odd} H_{k}(\Z_{\OO}; \mathbb{C}). \]
(for the precise definition of the centralizer $\Z_{\OO}$ and the normalizer $\N_{\OO}$ of $\OO$, see \ref{ZNO}).
\end{st}

\emph{proof}: From \ref{redloops} and \ref{1elliptic} we have $HP_{*}(\conv)_{\OO}=HP_{*}(\N_{\OO}; \A^{\natural}_{(\OO)})$. From the description of $\A^{\natural}_{(\OO)}$ (see \ref{ZNO}) it is not difficult to see that $\A^{\natural}_{(\OO)}= (\C_{V}^{\infty})^{\natural}_{|\OO}$ where $V$ is an open neighborhood of $\OO$ in $\nB{0}$ which is a submanifold of $\G$ (the existence of $V$ is ensured by the first part of \ref{stability}). Using \ref{stability} for $\f= id$ and the fact that $0 \rmap \mathbb{C} \rmap \Omega^{0}_{M}|_{A} \rmap \Omega^{0}_{M}|_{A} \rmap . \ . \ . \ $ is a c-soft resolution of $\mathbb{C} \in Sh(A)$ for any manifold $M$ and any $A \subset M$ (\cite{Iv}) we get, as in the previous proof, the relations expressing $HP_{*}$ in terms of the homology of $\N_{\OO}$. The passage to $\Z_{\OO}$ is ensured by the spectral sequence for the projection map $\f : \Z_{\OO} \rmap \N_{\OO}$ (it degenerates since $\ga/\f \simeq \ <\ga>$ is a finite cyclic group for all $\ga \in \OO$).  
\\ \newline
\par Remark that \ref{ultimaa} and this result ensure us that, for any $\OO$ as in \ref{last}, such that $\OO$ is closed in $\nB{0}$, there is a localized Chern-Connes character:
\[ Ch_{\OO}: K_{*}^{alg}(\conv) \rmap H_{*}(\Z_{\OO}; \mathbb{C}) .\]
In important situations it becomes an isomorphism after tensoring with $\mathbb{C}$; this, the problem of defining $Ch_{\OO}$ directly and explicitly, and also its usefulness in index problems will be considered elsewhere.

\begin{num}
{\bf Definition}: We call $\ga \in \nB{0}$ stable if its germ  $\sigma_{\ga}; (\nG{0}, s\ps \ga \pd ) \rmap (\nG{0}, s\ps \ga \pd)$ (as defined in \ref{sigma}) is stable (see \ref{stability}). We call $\OO \subset \nB{0}$ stable if it is \G-invariant and every $\ga \in \OO$ is stable.
\end{num}

\begin{st}
\label{lastt}
(hyperbolic case): For any smooth \'etale groupoid \G, and any hyperbolic $\OO \subset \nB{0}$, $HC_{*}(\conv)_{\OO}$ is a module over $H^{*}(\N_{\OO};\mathbb{C})$ and the $S$ in the $SBI$-sequence is the (cap-) product by the Euler class $e_{\OO} \in H^{2}(\N_{\OO})$.
\par Moreover, if $\OO$ is stable:
\[ HH_{n}(\conv)_{\OO}= \bigoplus_{p+q=n} H_{p}(\Z_{\OO}; \Omega^{q}_{\nB{0}} |_{\OO}), \]
\[ HC_{n}(\conv)_{\OO}= \bigoplus_{p+q=n} H_{p}(\N_{\OO}; \Omega^{q}_{\nB{0}} |_{\OO}), \]
and the $SBI$-sequence is the Gysin sequence for $\Z_{\OO} \rmap \N_{\OO}$ (see \ref{gysin}).
\end{st}
\par Remark that if $\OO$ is stable, there is a neighborhood $V$ of
 $\OO$ in $\nB{0}$ which is a submanifold of $\G$. In particular
 $\Omega^{*}_{V} \in Sh(V)$ makes sense. The notation
 $\Omega^{*}_{\nB{0}} |_{\OO} $ should be understood as the
 restriction of $\Omega^{*}_{V}$ to $\OO$; when $\OO$ is open in
 $\nB{0}$, this is simply $\Omega^{*}_{\OO}$. 
 \par \emph{proof}: The first part is a consequence of \ref{redloops}, \ref{1hyperbolic}. We are left with the last two equalities. We prove the second one (the first can be proved in the same way or as a
 consequence of \ref{stability}, \ref{redloops} , \ref{redcross}). The
 proof is an improvement of the second spectral sequence in \ref{chyperbolic} for the cyclic category $(\Z_{\OO}, \thp)$. There we used the spectral sequence of $\pi^{\, '}_{2}:\Lambda_{\infty}
 \wedge \Z_{\OO} \rmap \N_{\OO}$ and the strong deformation retract $\ga / \pi^{\ '}_{2} \tilde{\leftarrow} \Lambda_{\infty} \wedge <\ga> = \Lambda_{\infty}$ for all $\ga \in \OO$. Assume for simplicity that $\OO$ is open in $\nB{0}$. We have from \ref{localization}:
\[ HC_{*}(\conv)_{\OO} = H_{*}(\Lambda_{\infty} \wedge \Z_{\OO};\A^{\natural}_{\OO}) .\]
Due to the $SBI$ trick (see \ref{important}) and \ref{stability}, we can replace $\A^{\natural}_{\OO}$ by the $\thp$-cyclic $\Z_{\OO}$-sheaf $\B^{\natural}=(\C^{\infty}_{\OO})^{\natural}$. This is in fact a
 cyclic $\N_{\OO}$-sheaf. We use (see \ref{HS}) that the spectral sequence for $\pi^{\, '}_{2}$ comes from an equality:
 \[ H_{*}(\Lambda_{\infty} \wedge \Z_{\OO}; \B^{\natural})=\mathbb{H}_{*}(\N_{\OO}; (\el \pi^{\, '}_{2})_{\, !}\B^{\natural}) ,\]
 where $\tilde{H}_{q}((\el \pi^{\, '}_{2})_{\,!}\B^{\natural}) \in Sh(\N_{\OO})$ has the stalk at $\ga \in \OO$: 
\[ H_{q}(\ga/\pi^{\,'}_{2}; \B^{\natural})=H_{q}(\Lambda_{\infty}; \B^{\natural}_{\ga})\] 
 (the last equality follows from \ref{inv}). We get in this way a quasi-isomorphism of complexes of sheaves on $\OO$: $(\el \pi^{\, '}_{2})_{\,!}\B^{\natural} \simeq (\B^{\natural}, b)$, $b=$the
 Hochschild boundary. This is in $Sh(\N_{\OO})$ so:
 \[ HC_{*}(\conv)_{\OO}= \mathbb{H}_{*}(\N_{\OO}; (\B^{\natural}, b)) .\]
 Using once more \ref{stability} (for $\al=id$):
 $(\B^{\natural}, b) \simeq (\Omega_{\OO}^{*}, 0)$ and the second
 equality follows.  
\begin{num}
{\bf Corollary}(compare to \cite{Ni1}): If $\OO$ is hyperbolic and $\N_{\OO}$ has finite cohomological dimension, then $HP_{*}(\conv)_{\OO}=0$.
\end{num}


\subsection{The case of cohomology}

We need some facts about sheaves in order to make possible an analogous treatment of the cyclic \emph{co}homology for smooth \'etale groupoids.

\begin{num}
\label{cohomol}
\emph{
{\bf Cohomology via Bar-complexes}: We will need a version of \ref{barcomp}, \ref{homol} for computing the cohomology $H^{*}(\G; \A)$ for an \'etale groupoid \G\ and a left \G-sheaf \A\ (i.e.  $\A \in Sh(\G^{op})$ cf. \ref{sheaves}). This is well known (see \cite{SGA,Mo4} or \cite{Ha3} for a direct approach); we briefly review it. By definition, $H^{*}(\G; -)$ are the right derived functors of $\Gamma^{\, \G}=Hom_{Sh(\G^{op})}(\mathbb{C}, -)$.
We have a sequence of \'etale groupoids:
\[ \nG{q} \stackrel{\alpha_{q}}{\rmap}  \nG{0}  \stackrel{i}{\rmap} \G , \ \ \alpha_{q}(g_{1}, . . . , g_{q})=t(g_{1}), \ \ i={\rm the \ \  inclusion} \]
(here $\nG{q}, \nG{0}$ are spaces, viewed as \'etale groupoids as in \ref{exgr}.1). We have the induced functors:
\[ \xymatrix {
 Sh(\nG{q})  \ar@<1ex>[r]^{\alpha_{q *}} & Sh(\nG{0}) \ar@<1ex>[r]^{i_{*}} \ar@<1ex>[l]^{\alpha_{q}^{*}} & Sh(\G^{op}) \ar@<1ex>[l]^{i^{*}} } \]
and $i_{*}$ is right adjoint to $i^{*}$. Denoting $\C =i_{*}i^{*}: Sh(\G^{op}) \rmap Sh(\G^{op})$ and using the  unit of the adjunction $\eta: Id \rmap \C$ we get for any $\A \in Sh(\G^{op})$, a resolution of $\A$ by (left) \G-sheaves:
\[ \xymatrix {
\A \ar[r] & \C(\A) \ar@<-1ex>[r] \ar[r] & \C^{2}(\A) \ar@<-1ex>[r] \ar[r] \ar@<1ex>[r] & \C^{3}(\A) \ . \ . \ .  } \ \ .\]
From the spectral sequences induced by $i_{*}, \alpha_{q *}$ we have for any $\A \in Sh(\G^{op}), \B \in Sh(\nG{0})$:
\[ H^{*}(\G; i_{*}i^{*}\A)=H^{*}(\nG{0}; i^{*}\A)\ \ , \ \ H^{*}(\nG{0}; \alpha_{q *}\alpha_{q}^{*}\B)=H^{*}(\nG{q}; \alpha_{q}^{*}\B) .\]
Using also that $(i^{*}i_{*})^{q}i^{*}=\alpha_{q *}\alpha_{q}^{*}i^{*}$, we get:
\[ H^{*}(\G; \C^{q}(\A))=H^{*}(\nG{q-1}; \alpha_{q-1}^{*}\A) .\]
We say that $\A$ is $\G$-acyclic if the left side of the previous equality is $0$ for each $*>0, q \geq 1$; in this case we get that $\A \rmap \C^{\dt}(\A)$ is a resolution of $\A$ by $\Gamma^{\, \G}$-acyclic objects, so it can be used for computing $H^{*}(\G; \A)$. In other words, if $\A$ is $\G$-acyclic, then $H^{*}(\G; \A)$ is computed by a co-simplicial vector space with 
\[ B^{\dt}(\G; \A)=\Gamma^{\G}(\C^{\dt+1}(\A))=H^{0}(\G; \C^{\dt+1}(\A))=\Gamma(\nG{\dt}; \alpha_{\dt}^{*}\A) \] 
(after spelling out the co-boundaries, we get  formulas which are similar to the ones in \ref{barcomp}). More generally, if $\es \in Sh(\G^{op})$, then $H^{*}(\G; \es)$ is computed by the the double complex $B^{\dt}(\G; \A^{\dt})$ where $\es \rmap \A^{\dt}$ is any resolution in $Sh(\G^{op})$ by $\G$-acyclic sheaves. The same disscusion carries over to \'etale categories. 
}
\end{num}

\begin{num}
\label{dual}
\emph{
{\bf Dual sheaves}: Let $M$ be a space, $\A \in Sh(M)$ c-soft. The correspondence $U \mapsto \gc(U;\A)^{\vee}$ (recall that "$\vee$" stands for the algebraic dual) defines a sheaf on $M$, denoted $\A^{\vee}$; it is flabby (even injective; see \cite{Iv}). This is well known in the Hausdorff case and carries over to the general case (\cite{CrMo}). The construction has the property that for any $f:M \rmap N$ continuous, and for any c-soft sheaf $\A \in Sh(M)$:
\[ f_{\ !}(\A)^{\vee}=f_{*}(\A^{\vee}) .\]
When applied to sheaves on an \'etale category, it gives a correspondence:
\[ Sh(\G) \ni \A  \ ,\ c-{\rm soft} \ \mapsto \A^{\vee} \in Sh(\G^{op}) \ ,\   \G-{\rm acyclic} ,\]
and $B_{\dt}(\G; \A)^{\vee}= B^{\dt}(\G; \A^{\vee})$.
}
\end{num}

\begin{num}
\label{currents}
\emph{
{\bf Currents}: Let $M$ be a manifold (not necessarily Hausdorff) of dimension $n$. There is an obvious notion of $q$-currents: $\Omega_{q}(M)=\{ u: \Omega_{c}^{q}(M) \rmap \mathbb{C}: u\, |_{U} \in \Omega_{q}(U) \ \ $ for all local coordinates charts $\ U \}$ ($\Omega_{q}(U)$ has the usual meaning if $U$ is Hausdorff). We get the sheaf of $q$-currents on $M$ which is a kind of topological dual of the sheaf $\Omega_{M}^{q}$; we denote it by $(\Omega_{M}^{q})^{'}$. In general we have two different resolutions of the (complex) orientation sheaf: $or_{M} \rmap (\Omega_{M}^{n})^{'} \rmap  (\Omega_{M}^{n-1})^{'} \rmap . . . $ and $or_{M} \rmap (\Omega_{M}^{n})^{\vee} \rmap  (\Omega_{M}^{n-1})^{\vee} \rmap . . . $. (\cite{BoTu,Iv}) and an obvious "forgetting continuity" map of complexes in $Sh(M)$:
\[  ((\Omega_{M}^{*})^{\, '} , d_{dRh}^{\, '} ) \rmap ((\Omega_{M}^{*})^{\vee} , d_{dRh}^{\, \vee} ) .\]
Applying to this the global sections functor, we get two cohomologies: de Rham cohomology with coefficients in $or_{M}$ (or, equivalently, closed de Rham homology) and the sheaf cohomology of $or_{M}$, together with a linear map between them:
\[ H^{*}_{dRh}(M; or_{M}) \rmap H^{*}(M; or_{M}) .\]
In the Hausdorff case, all sheaves $(\Omega_{M})^{'}$ and $(\Omega_{M})^{\vee}$ are soft and the map above is an isomorphism; in the non-Hausdorff case, $(\Omega_{M})^{'}$ may fail to be acyclic for cohomology and the two cohomologies are not isomorphic in general. The same discussion carries over to the case of \'etale groupoids (see \cite{Ha1}).
}
\end{num}

\begin{num}
{\bf Definition}: Let \G\ be a smooth \'etale groupoid. Define $HH^{*}(\conv), HC^{*}(\conv)$,\\ $HP^{*}(\conv)$ as the cohomologies of $(\conv^{\natural})^{\vee}$.
\end{num}
\ \ \ Emphasize that these are not the algebraic cohomologies of $\conv$ since in the definition of $\conv^{\natural}$ we used the topology.
 It is possible also to use $(\conv^{\natural})^{'}$ for defining the cyclic cohomologies. Because of the facts explained in \ref{currents} they give the same \emph{periodic} cyclic cohomologies in the case of Hausdorff groupoids. In the general case the difference between them is the same as the difference between the two cohomologies described in \ref{currents}. 
\par Remark that we chose the maximal definition such that we keep the pairing with the cyclic homology and such that it is a receptacle for "Chern-Connes character" maps (\cite{Co3,Ni2}). And, as we shall see, it is "computable".

\begin{st}
\label{cohelliptic}
: For any smooth \'etale groupoid \G, and any elliptic open $\OO \subset B^{\ps 0 \pd}$ which is a topological manifold (not necessarily Hausdorff) of dimension $q$:
\[ HP^{0}(\conv)_{\OO}=\bigoplus_{k-even} H^{k+q}(\N_{\OO}; or_{\OO}),\]
\[ HP^{1}(\conv)_{\OO}=\bigoplus_{k-odd} H^{k+q}(\N_{\OO}; or_{\OO} ), \]
and the pairing with $HP_{*}(\conv)_{\OO}$ is in fact the Poincar\'e-duality pairing (see \ref{poincare}).
\par In particular, the Connes pairing $HP^{*} \times HP_{*} \rmap \mathbb{C}$ induces an inclusion:
\[ HP^{*}(\conv)_{\OO} \hookrightarrow HP_{*}(\conv)_{\OO}^{\vee} .\]
\end{st}

\emph{proof}: We explain how to make the analogy with the case of homology (and get more of the lemmas we need as consequences of that case). The analogue of \ref{redloops} we get for free (as a consequence). The analogy with the sub-paragraph 2.4 is ensured by \ref{cohomol} and \ref{dual}. To see this, keep the notations from \ref{cycrpr}, \ref{NNN}. Since our $\A=\C_{\nG{0}}^{\infty}$ is c-soft and  
since $(\A \cross \G )^{\natural}(n)$ is isomorphic to $\gc (\Z^{\ps n \pd}; \A^{\ps n \pd}(\nB{0}))$ (see the isomorphism in the proof of \ref{redcross}) we get that $((\A \cross \G )^{\natural}(n))^{\vee}$ is isomorphic to $\gc (\Z^{\ps n \pd}; \A^{\ps n \pd}(\nB{0})^{\vee})$ \emph{and} each $\A^{\ps n \pd}(\nB{0})^{\vee}$ is injective (cf. \ref{dual}). From this point on, the analogy is ensured by \ref{cohomol} and the isomorphism in \ref{dual}. The analogue of \ref{stability} we also get for free (just by dualizing) and we get down to the resolutions of the orientation sheaf given by the duals of the sheaves of forms. This also introduces the shifting in the result. The identification with the Poincar\'e pairing is straightforward.\\
\par We leave the statement of the hyperbolic case to the reader. As an obvious (but interesting) consequence we have:

\begin{num}
\label{of}
{\bf Corollary}: Let $\OO \subset \nB{0}$ be \G-invariant such that it is a disjoint union of manifolds. Assume that $\N_{\OO}$ has finite cohomological dimension.Then:
\begin{enumerate}
\item if $\OO$ is elliptic: $HP^{*}(\conv)_{\OO} \simeq HP_{*}(\conv)_{\OO}^{\vee} $ (induced by the Connes pairing),
\item if $\OO$ is hyperbolic: $HP^{*}(\conv)_{\OO} = HP_{*}(\conv)_{\OO} = 0$.
\end{enumerate}
\end{num}

\begin{num}
\emph{
{\bf At units}: As a particular case of \ref{cohelliptic} we have the localization at units:
\[ HP^{*}(\conv)_{\pss 1 \pdd}=\bigoplus_{k\equiv * (mod 2)} H^{k+q}(\G; or ) .\] 
By using \ref{classif}, we get the connection with the cohomology of the classifying space. But remark that the result in this form, as the cohomology of \G, gives us more freedom in constructing cocycles. For instance, any \G-vector bundle $\xi$ gives its Chern classes $c_{*}(\xi) \in H^{*}(\G)$ in explicit cocycles (we can repeat the construction of Chern classes as done in \cite{Bo}).
}
\end{num}

\subsection{Group Actions on Manifolds}

We use the previous results to describe the homologies of the cross-product (locally convex) algebra $C_{c}^{\infty}(M) \cross G$. Here $G$ is a discrete group acting smoothly on the Hausdorff manifold $M$. This algebra coincides with the convolution algebra of the smooth \'etale groupoid $\G=M \cross G$ (see \ref{exgr}). For any $g \in G$, denote by $Z_{g}=\{ h\in G: hg=gh \}$ and $N_{g}=Z_{g}/<g>$ the centralizer and the normalizer of $g$, and by $M^{g}$ the points fixed by $g$. Denote by "$\sim$" the conjugacy relation on $G$ and put $<G> =G/\sim$. The loop space of \G\ is $\nB{0}=\{ \ps x, g\pd \in M \times G : xg=x \}$ and is usually denoted by $\widehat{M}$ (\cite{BaCo}). Any $g \in G$ defines an invariant open $\OO_{g}=\{ (x, h) \in \widehat{M} : h \sim g \}$ and $\widehat{M}=\coprod_{g \in <G>} \OO_{g}$. In particular we have the well-known decomposition (see \cite{Bry}, \cite{FeTs}, \cite{Ni1}):
\[ HH_{*}( C_{c}^{\infty}(M) \cross G)= \bigoplus_{g \in <G>} HH_{*}( C_{c}^{\infty}(M) \cross G)_{\ps g \pd},  \]
and the analogues for $HC_{*}, HP_{*}$. For any $g$ we have obvious Morita equivalences $\Z_{\OO_{g}} \simeq M^{g} \cross Z_{g}$,  $\N_{\OO_{g}} \simeq M^{g} \cross N_{g}$.
\par  In the elliptic case we get:
\[ HP_{*}( C_{c}^{\infty}(M) \cross G)_{\ps g \pd}=\prod_{k \equiv * (mod2)} H_{k}(M^{g} \cross N_{g}; \mathbb{C}) ,\]
which, with \ref{poincare} and \ref{classif} in mind, is similar to the description given in \cite{BrNi}. Remark also that if $G$ acts properly on $M$ we get (using \ref{HS} for $M^{g} \cross N_{g} \rmap M^{g}/N_{g}$ and \ref{lili}):
\[ HP_{*}(\C_{c}^{\infty}(M) \cross G)= \bigoplus_{g \in <G>} \prod_{k=* mod 2} H_{c}^{k}(M^{g}/Z_{g}) \]
as the target of the Chern character. This agrees with the definition given in \cite{BaCo}. 
\par Assume now that $g \in G$ is hyperbolic. We have a projection $(\Z_{\OO_{g}}, \thp) \rmap (Z_{g}, g)$ of hyperbolic cyclic groupoids (see example 3 in \ref{excyclic}) which induces a map:
\[ H^{*}_{N_{g}}(pt; \mathbb{C})=H^{*}(N_{g}; \mathbb{C}) \rmap H^{*}_{N_{g}}(M^{g}; \mathbb{C})=H^{*}(\N_{\OO_{g}}; \mathbb{C}) .\]
From \ref{gysin}, $e_{\OO_{g}} \in H^{2}_{N_{g}}(M^{g}; \mathbb{C})$ is the image by this map of the Euler class of $(Z_{g}, g)$. This Euler class is denoted by $e_{g} \in H^{2}(N_{g})$; by the previous work of Nistor (\cite{Ni1}) we know it is represented by the extension $<g> \ \rmap Z_{g} \rmap N_{g}$. We get the following :

\begin{num}
\label{off}
{\bf Corollary}: For $g \in G $ hyperbolic, $HC_{*}(\C_{c}^{\infty} \cross \G)_{\ps g \pd}$ is a module over the Borel cohomology ring $H^{*}_{N_{g}}(M^{g}; \mathbb{C})$ and $S$ in the $SBI$-sequence is the product by the image of the Euler class $e_{g} \in H^{2}(N_{g}; \mathbb{C})$ in $H^{2}_{N_{g}}(M^{g}; \mathbb{C})$.
\par Moreover, if $g$ acts stably (for instance if it preserves a metric) then:
\[ HH_{n}(\C_{c}^{\infty} \cross \G)_{\ps g \pd}= \bigoplus_{p+q=n} H_{p}(Z_{g}; \Omega^{q}_{c}(M^{g})) ,\]
\[ HC_{n}(\C_{c}^{\infty} \cross \G)_{\ps g \pd}= \bigoplus_{p+q=n} H_{p}(N_{g}; \Omega^{q}_{c}(M^{g})) ,\]
and the $SBI$-sequence is the sum of the Gysin-sequences for $<g> \ \rmap Z_{g} \rmap N_{g} $.
\end{num}

\ \ \ Compare to \cite{Bu,Ni1}. Also we have the dual results for cohomology. In particular:

\begin{num}
{\bf Corollary}: If $g \in G$ is hyperbolic and the image of $e_{g} \in H^{2}(N_{g}; \mathbb{C})$ in $H^{2}_{N_{g}}(M^{g}; \mathbb{C})$ is nilpotent, then $HP_{*}( C_{c}^{\infty}(M) \cross G)_{\ps g \pd}$ and $HP^{*}( C_{c}^{\infty}(M) \cross G)_{\ps g \pd}$ vanish.
\end{num}


\subsection{Foliations}

\ \ \ Let $(M, \F)$ be a foliated manifold. The holonomy groupoid (see \ref{exgr}), $Hol(M, \F)$ is a smooth groupoid which is non-Hausdorff in general. When restricted to a complete transversal $T$ it becomes an \'etale groupoid $Hol_{T}(M, \F)$. As proved in \cite{HiSk}, the choice of $T$ is not important when we talk about the associated non-commutative spaces (i. e. the associated $C^{*}$-algebras). The case of their smooth convolution algebras seems to be slightly more difficult. However we have:

\begin{st}
\label{folinv}
: For any foliated manifold $(M, \F)$ the correspondences:
\[ T \mapsto HH_{*}, HC_{*}, HP_{*}, HH^{*}, HC^{*}, HP^{*}  ( \C_{c}^{\infty}( Hol_{T}(M, \F))) \]
do not depend on the choice of the complete transversal $T$; so they give well defined invariants of the (``leaf space'' of the) foliation. These are denoted by $HH_{*}(M, \F), HC_{*}(M, \F), HP_{*}(M, \F)$.
\end{st}

\ \ \ Remark that this is not implied by the results of \cite{Mr}, namely: in the \emph{Hausdorff} case, for any two complete transversals $T, T\,'$, there is an \emph{algebraic} Morita equivalence: $\C_{c}^{\infty}( Hol_{T}(M, \F)) \simeq \C_{c}^{\infty}( Hol_{T\,'}(M, \F))$.
\par
\par \emph{proof of \ref{folinv}}: Let $T, T\,'$ be two complete transversals; $\G=Hol_{T}(M, \F), \N$ and $\B=\A^{\natural}_{\nB{0}}$ be the constructions from \ref{NNN} for \G; the similar constructions for $T\,'$ are denoted by $\G\,', \N\,', \B\,'$. Replacing $T$ by $T \coprod T\,'$ if needed, we can assume $T \subset T\,'$. We have continuous functors $\G \rmap \G\,', \N \rmap \N\,'$ which are Morita equivalences. Since $HH_{*}(\conv)=\mathbb{H}_{*}(\N; (\B, b))$ (cf \ref{redloops}, \ref{redcross}, \ref{important}) and the analogue for $\G\,'$, it is enough to use Morita invariance for homology and an $SBI$ argument.  
\begin{num}
\emph{
{\bf At units}: For the localization at units we get:
\[ HP^{*}(M/\F)_{\ps 1 \pd}=\bigoplus_{k \equiv q+* mod 2} H^{k}(M/\F; or) .\]
This is a common point between two different approaches to model the leaf-space $M/\F$ as a generalized space: one in the spirit of non-commutative geometry and one in the spirit of Grothendieck, by looking at all the sheaves on $Hol(M; \F)$ (see \cite{Mo1}). In particular, the right side is the cohomology of the orientation sheaf of $M$ inside the category of sheaves of complex vector spaces on $Hol(M, \F)$. See also \cite{Mo3} for the connection with De Rham (or basic) cohomology of the leaf-space (\cite{Mol}).
}
\end{num}

\begin{num}
{\bf Examples}:
\end{num}

\begin{enumerate}
\item Consider the foliation of the open Moebius band and the complete transversal $T$ as in the picture. Then $\G=Hol_{T}(M, \F)$ has $\nG{0}=(-1, 1), \nG{1}=(-1, 1) \cup_{D} (-1, 1)$ where $D=(-1, 0) \cup_{} (0, 1)$. It is not Hausdorff and has just elliptic loops. Applying \ref{last} we get:
\[ HP_{0}(M/\F)=0\ , \ HP_{1}(M/\F)=\mathbb{C}\ .\]
\setlength\epsfxsize{9cm}
\vskip 0cm
\hskip 2.7cm
\epsffile{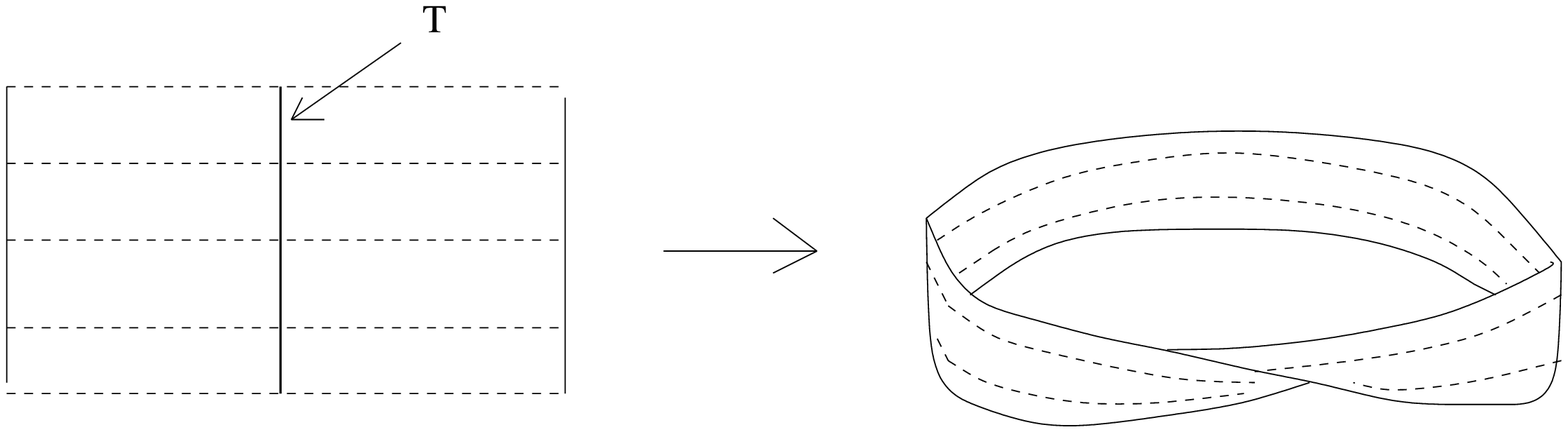}

\item Consider the Kronecker foliation on the torus (see I.4.$\beta$ in \cite{Co3}) which comes from the foliation of the plane by lines of slope $\al=2\pi \thp, \ \thp \in \mathbb{R}-\mathbb{Q}$ . Choosing $T$ as in the picture, the reduced groupoid becomes $S^{1}\cross \mathbb{Z}$ where $\mathbb{Z}$ acts on $S^{1}$ by rotations by $\al$. The only elliptic loops are the units. From \ref{last} we get $HP_{*}(M/\F)=HP_{*}(M/\F)_{\ps 1 \pd}$ and from \ref{ultima} this is computable in terms of $H_{*}(S^{1} \cross \mathbb{Z})$. An easy algebraic computation shows that the last group is $\mathbb{C}$ if $* \in \{ -1, 1 \}$, $\mathbb{C} \oplus \mathbb{C}$ if $*=0$ and $0$ otherwise. So:
\[ HP_{0}(M/\F)=HP_{1}(M/\F)= \mathbb{C} \oplus \mathbb{C} .\]
A similar result is obtained for cohomology. Compare to theorem 53 in \cite{Co2}.
\setlength\epsfxsize{9cm}
\vskip
3mm\hskip 2.7cm
\epsffile{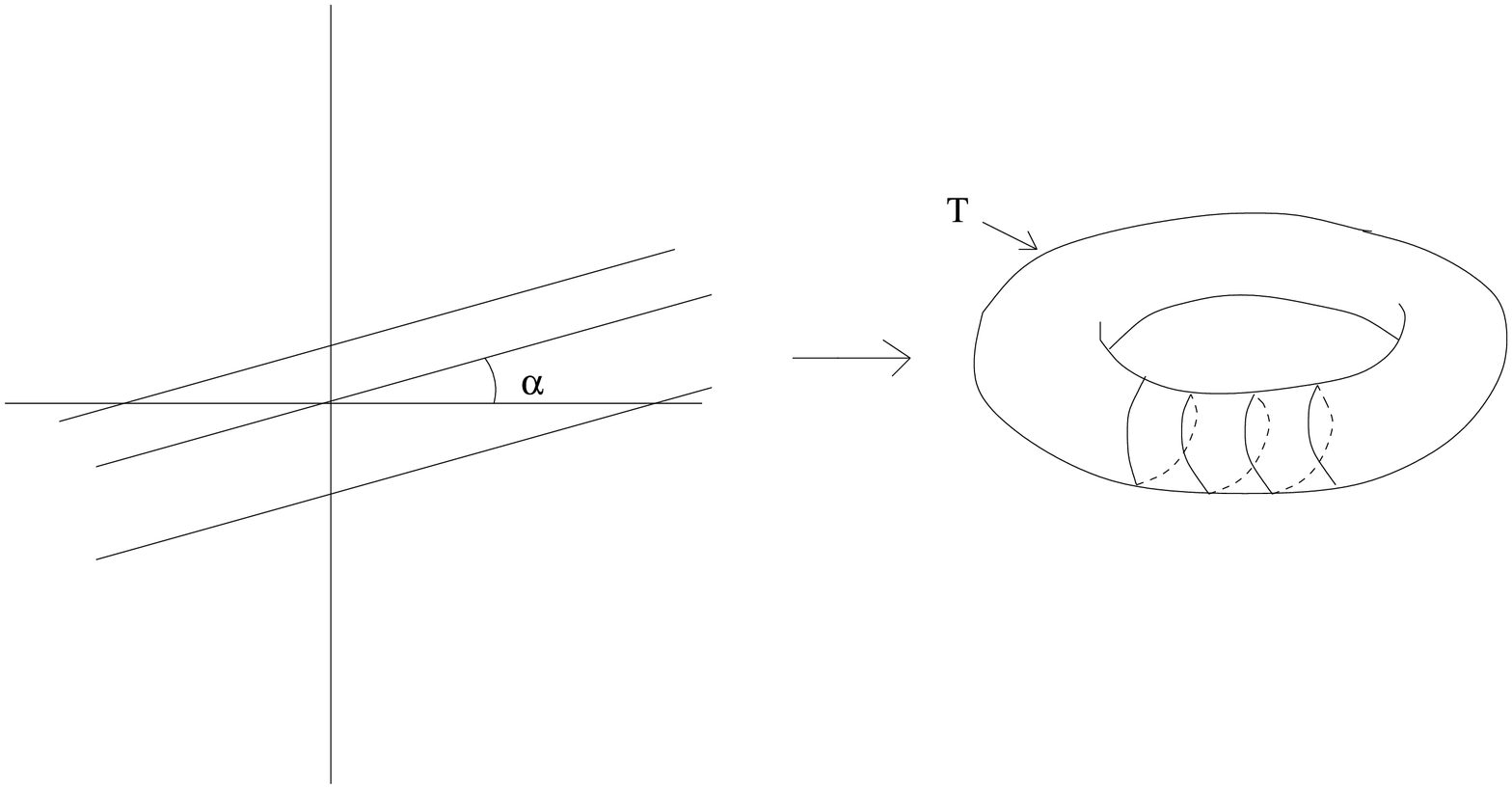}

\item For Haefliger's groupoid $\Gamma^{q}$ (which is non-Hausdorff), the localization at units gives the homology and cohomology of $\Gamma^{q}$. The last one gives universal cocycles which induce characteristic classes for foliations of codimension $q$. This is the case for the Godbillon-Vey invariant or, more generally, classes coming from Gelfand-Fuchs cohomology. See \cite{Ha1} and \cite{Co4}. We get in this way an interpretation in terms of cyclic cocycles for these cohomology classes. 
\end{enumerate}


\end{document}